%% file: Arikan_Binding_No3_Pr_GGT.tex
\documentclass[reqno]{gokova}
\usepackage{epsfig,graphicx}
\usepackage{amssymb}

\begin{document}
\input goksty.tex      

\setcounter{page}{90}
\volume{13}

\newcommand{\A}{$Aut(\Sigma,\partial \Sigma) \,\,$}
\newcommand{\D}{$Dehn^{+}(\Sigma,\partial \Sigma) \,\,$}
\newcommand{\V}{$Veer(\Sigma,\partial \Sigma) \,\,$}
\newcommand{\rv}{right-veering }

\newcommand{\au}{$Aut(\Sigma,\partial \Sigma)$}
\newcommand{\de}{$Dehn^{+}(\Sigma,\partial \Sigma)$}
\newcommand{\ve}{$Veer(\Sigma,\partial \Sigma)$}


\title[Planar Contact Structures with Binding Number Three]
{Planar Contact Structures with Binding Number Three}
\author[Mehmet F. Ar\i kan]{ Mehmet F\i rat Ar\i kan }

\thanks{The author was partially supported by NSF Grant DMS0244622}

\address{Department of Mathematics, Michigan State University,
          East Lansing, MI 48824-1027}
\email{arikanme@msu.edu}

\begin{abstract}
In this article, we find the complete list of all contact structures
(up to isotopy) on closed three-manifolds which are supported by an
open book decomposition having planar pages with three (but not
less) boundary components. We distinguish them by computing their
first Chern classes and three dimensional invariants (whenever
possible). Among these contact structures we also distinguish tight
ones from those which are overtwisted.
\end{abstract}

\maketitle

\section{Introduction} \label{introduction}
\noindent Let $(M,\xi)$ be a closed oriented 3-manifold with the contact
structure $\xi$, and let $(S,h)$ be an open book (decomposition) of
$M$ which is compatible with $\xi$. In this case, we also say that
$(S,h)$ supports $\xi$ (for the definitions of these terms see the
next section). Based on Giroux's correspondence theorem (Theorem
\ref{Giroux}), two natural questions have been asked in \cite{EO}:

\medskip \noindent (1) What is the possible minimal page genus $g(S)=$genus$(S)$?

\medskip \noindent (2) What is the possible minimal number of boundary components of a
page $S$ with $g(S)$ minimal?

\medskip \noindent In \cite{EO}, two topological invariants $sg(\xi)$ and $bn(\xi)$
were defined to be the answers. More precisely, we have:
\[
sg(\xi)=\min \{ \, \,g(S) \, \, \vert \, \,(S,h) \text{ an open book
decomposition supporting } \xi\},
\]
called \emph{the support genus} of $\xi$, and
\[
bn(\xi)=\min \{ \, \, |\partial S|  \, \,\vert  \, \, (S,h) \text{
an open book decomposition supporting }  \xi \text{ and }
g(S)=sg(\xi)\},
\]
called \emph{the binding number} of $\xi$. There are some partial results
about these invariants. For instance, it is proved in \cite{Et1} that
if $(M,\xi)$ is overtwisted, then $sg(\xi)=0.$

\medskip \noindent Unlike the overtwisted case, there is not much known yet
for $sg(\xi)$ if $\xi$ is tight. The algorithm given in \cite{Ar}
finds a reasonable upper bound for $sg(\xi)$  using the given
contact surgery diagram of $\xi$. However, there is no systematic
way to obtain actual values of $sg(\xi)$ and $bn(\xi)$ yet.

\medskip \noindent One of the ways to work on the above
questions is to get a complete list of contact manifolds
corresponding to a fixed support genus and a fixed binding number.
To get such complete list, we consider all possible monodromy maps
$h$. The first step in this direction is the following result given
in \cite{EO}. Throughout the paper $L(m,n)$ stands for the lens
space obtained by $-m/n$ rational surgery on an unknot.
\begin{thm} [\cite{EO}] \label{bindingnumber=<2}
Suppose $\xi$ is a contact structure on a 3-manifold $M$ that is
supported by a planar open book (i.e., $sg(\xi)=0$). Then
\begin{enumerate}
\item If $bn (\xi)=1$, then $\xi$ is the standard tight contact structure on $S^3.$
\item If $bn (\xi)=2$ and $\xi$ is tight, then $\xi$ is the unique tight contact structure on the lens space $L(m,m-1)=L(m,-1)$ for some $ m \in {\Z_+ \cup \{0\}}$.
\item If $bn(\xi)=2$ and $\xi$ is overtwisted then $\xi$ is the overtwisted contact structure on $L(m,
1)$, for some $m \in {\Z_+}$, with $e(\xi)=0$ and $d_3(\xi)=-\frac
14 m +\frac 34$ where $e(\xi)$ and $d_3(\xi)$ denotes the Euler
class and $d_3-$invariant of $\xi$, respectively. When $m$ is even
then the refinement of $e(\xi)$ is given by
$\Gamma(\xi)(\mathfrak{s})=\frac m2$ where $\mathfrak{s}$ is the
unique spin structure on $L(m,1)$ that extends over a two handle
attached to a $\mu$ with framing zero. Here we are thinking of
$L(m,1)$ as $-m$ surgery on an unknot and $\mu$ is the meridian to
the unknot.
\end{enumerate}
\end{thm}
\noindent We remark that Theorem \ref{bindingnumber=<2} gives the
complete list of all contact $3$-manifolds which can be supported by
planar open books whose pages have at most 2 boundary components.
Next step in this direction should be to find all contact
$3$-manifolds $(M,\xi)$ such that $sg(\xi)=0$ and $bn (\xi)=3$. In
the present paper, we will get all such contact structures, and also
distinguish tight ones by looking at the monodromy maps of their
corresponding open books (See Theorem \ref{Contact_Type} and Theorem
\ref{the_List}). After the preliminary section (Section
\ref{preliminaries}), we prove the main results in Section
\ref{Results}. Although some ideas in the present paper have been
already given or mentioned in \cite{EO}, we will give their explicit
statements and proofs in our settings. We finish
this section by stating the main results.

\begin{figure}[ht]
  \begin{center}
   \includegraphics{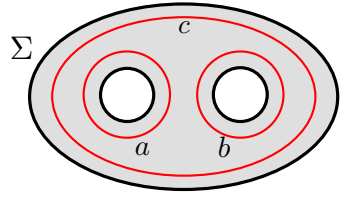}
   \caption{The surface $\Sigma$ and the curves giving the generators of \au.}
  \label{our surface}
    \end{center}
\end{figure}

\medskip \noindent Let $\Sigma$ be the compact oriented surface with $|\partial
\Sigma|=3$, and consider the boundary parallel curves $a,b,c$ in
$\Sigma$ as in the Figure \ref{our surface}. Through out the paper,
$\Sigma$ will always stand for this surface whose abstract picture
is given below. Let \A be the group of (isotopy classes of)
diffeomorphisms of $\Sigma$ which restrict to the identity on
$\partial \Sigma$. (Such
diffeomorphisms are automatically orientation-preserving). \\

\noindent It is known (see \cite{Bi}) that
\begin{center}
 \A$=\Z \langle D_a \rangle \oplus \,\,\Z \langle D_b \rangle \oplus \,\,\Z \langle D_c \rangle \cong {\Z}^3$
\end{center}
where $D_a, D_b, D_c$ denote positive Dehn twists along the curves
$a,b,c$ given as in Figure \ref{our surface}. In the rest of the paper, we will not make any distinction between
isotopy classes of arcs/curves/maps and the individual
arcs/curves/maps.

\medskip \noindent We start with studying the group \A in details.
Since generators commute with each other, we have that
\medskip
\begin{center}
 \A$=\{{D_a}^p{D_b}^q{D_c}^r | \, p,q,r \in \Z\}$.
\end{center}

\begin{figure}[ht]
  \begin{center}
   \includegraphics{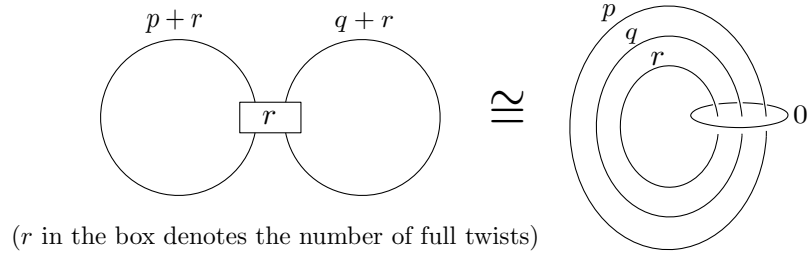}
   \caption{Seifert fibered manifold $Y(p,q,r)$.}
  \label{our model}
    \end{center}
\end{figure}

\medskip \noindent
For any given $p,q,r \in \Z$, let $Y(p,q,r)$ denote the smooth
$3$-manifold given by the smooth surgery diagram in Figure \ref{our
model} (diagram on the left). It is an easy exercise to check that
$Y(p,q,r)$ is indeed diffeomorphic to Seifert fibered manifold
given in Figure \ref{our model} (diagram on the right).

\noindent Now we state the following theorem characterizing all
closed contact $3$-manifolds whose contact structures supported by
open books $(\Sigma,\phi={D_a}^p{D_b}^q{D_c}^r)$.

\begin{figure}[ht]
  \begin{center}
   \includegraphics{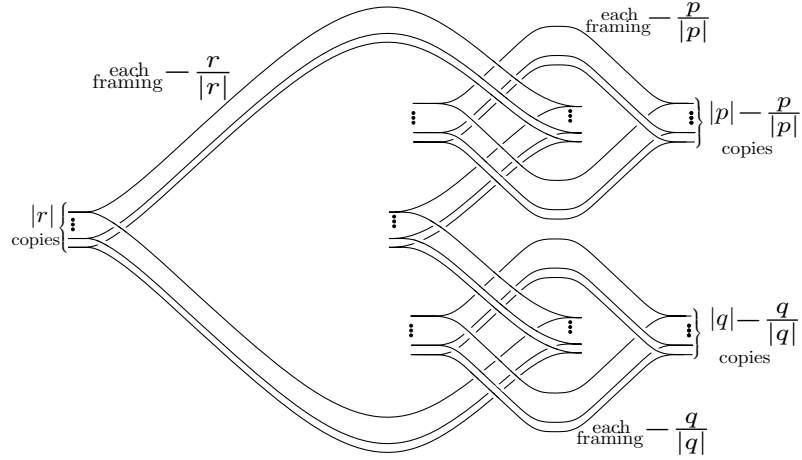}
   \caption{Contact manifold $(Y(p,q,r),\xi_{p,q,r})$.}
  \label{modelstructure}
    \end{center}
\end{figure}

\begin{thm}  \label{Contact_Type}
Let $(M,\xi)$ be a contact manifold supported by the open book
$(\Sigma,\phi)$ where $\phi={D_a}^p{D_b}^q{D_c}^r \in$ \A for fixed
integers $p,q,r$. Then $(M,\xi)$ is contactomorphic to
$(Y(p,q,r),\xi_{p,q,r})$ where $\xi_{p,q,r}$ is the contact
structure on $Y(p,q,r)$ given by the contact surgery diagram in
Figure \ref{modelstructure}. Moreover,
\begin{enumerate}
\item $\xi$ is tight (in fact holomorphically fillable) if $p\geq0, q\geq0,r\geq0$, and
\item $\xi$ is overtwisted otherwise.
\end{enumerate}
\end{thm}


\begin{rem} \label{howtochange1}
In Figure \ref{modelstructure}, if $r=0$, then we completely delete
the family corresponding to $r$ from the diagram, so  we are left with two families of Legendrian
curves which do not link to each other, and so the contact surgery
diagram gives a contact structure on the connected sum of two lens
spaces. However, if $p=0$ (or $q=0$), then we replace the Legendrian
family corresponding to $p$ (or $q$) by a single Legendrian unknot
with $tb$ number equal to $-1$, and we do $(+1)$-contact surgery on
the new unknot. Note also that Figure \ref{modelstructure} is
symmetric with respect to $p$ and $q$. This reduces the number of
cases in the proof of Theorem \ref{the_List}.
\end{rem}

\medskip \noindent
Of course not all $\xi_{p,q,r}$ have binding number three:

\begin{thm}  \label{the_List}
Let $(M,\xi)$ be a closed contact $3$-manifold with $sg(\xi)=0$ and
$bn(\xi)=3$. Then $(M,\xi)$ is contactomorphic to some
$(Y(p,q,r),\xi_{p,q,r})$ satisfying the following conditions:
\begin{center}
\begin{enumerate}

\item If $\; r=0$, then $p\neq1$ and $q\neq1$.
\item If $\; r=1$, then $p\notin \{-1,0\}$ and $q\notin \{-1,0\}$.
\item If $\; r=-1$, then $p\neq1$ and $q\neq1$.
\item If $\, |\,r|\geq2$, then $p\,q\neq-1$ and $(p,q)\notin \{(1,0),(0,1)\}$.
\end{enumerate}

\end{center}
\end{thm}

\noindent
Suppose that $(M,\xi)$ is a closed contact $3$-manifold with $sg(\xi)=0$ and
$bn(\xi)=3$, and let $c_1=c_1(\xi) \in H^2(M;\Z)$ denote the first Chern
class, and $d_3=d_3(\xi)$ denote the $3$-dimensional invariant
$($which lies in $ \mathbb{Q}$ whenever $c_1$ is a torsion class in
$H^2(M;\Z))$. Using $c_1$ and $d_3$, we can distinguish these structures in most of the cases.
In fact, we have either $M$ is a lens space, or a connected sum of lens spaces,
or a Seifert fibered manifold with three singular fibers. If one of the first two holds,
then using the tables given in Section \ref{Results} and \ref{last section},
one can get the complete list of all possible $(M,\xi)$ without any repetition.
That is, the contact structures in the list are all distinct pairwise and
unique up to isotopy. On the other hand, if the third holds, we
can also study them whenever $c_1$ is a torsion class. More discussion
will be given in Section \ref{last section}.

\medskip \noindent {\em Acknowledgments.\/} The author would like to
thank Selman Akbulut, \c{C}a\u{g}r{\i} Karakurt, Burak \"Ozba\u{g}c\i,
and Andr\'as Stipsicz  for helpful conversations.

\section{Preliminaries} \label{preliminaries}
\subsection{Contact structures and open book decompositions}

\noindent A $1$-form $\alpha \in \Omega^1(M)$ on a $3$-dimensional
oriented manifold $M$ is called a \emph{contact form} if it satisfies
$\alpha \wedge d\alpha \neq 0$. An \emph{oriented contact structure}
on $M$ is then a hyperplane field $\xi$ which can be globally
written as kernel of a contact $1$-form $\alpha$. We will always
assume that $\xi $ is a \emph{positive} contact structure, that is,
$\alpha \wedge d\alpha > 0$. Two contact structures $\xi_0, \xi_1$
on a $3$-manifold are said to be \emph{isotopic} if there exists a
1-parameter family $\xi_t$ ($0\leq t\leq 1$) of contact structures
joining them. We say that two contact $3$-manifolds $(M_1,\xi_1)$
and $(M_2,\xi_2)$ are \emph{contactomorphic} if there exists a
diffeomorphism $f:M_1\longrightarrow M_2$ such that
$f_\ast(\xi_1)=\xi_2$. Note that isotopic contact structures give
contactomorphic contact manifolds by Gray's Theorem. Any contact
$3$-manifold is locally contactomorphic to $(\mathbb{R}^3,\xi_0)$
where \emph{standard contact structure} $\xi_0$ on $\mathbb{R}^3$
with coordinates $(x,y,z)$ is given as the kernel of
$\alpha_0=dz+xdy$. The standard contact structure $\xi_{st}$ on the
$3$-sphere $S^3=\{(r_1,r_2,\theta_1,\theta_2) : r_1^2+r_2^2=1\}
\subset \mathbb{C}^2$ is given as the kernel of
$\alpha_{st}=r_1^2d\theta_1+r_2^2d\theta_2$. One basic fact is that
$(\mathbb{R}^3, \xi_0)$ is contactomorphic to $(S^3 \setminus \{pt
\},\xi_{st})$. For more details on contact geometry, we refer the
reader to \cite{Ge}, \cite{Et3}.

\medskip \noindent
\noindent An \emph{open book decomposition} of a closed $3$-manifold
$M$ is a pair $(L, f)$ where $L$ is an oriented link in $M$, called
the \emph{binding}, and $f: M\setminus L \to S^1$ is a fibration such
that $f^{-1}(t)$ is the interior of a compact oriented surface $S_t
\subset M$ and $\partial \Sigma_t=L$ for all $t \in S^1$. The
surface $S=S_t$, for any $t$, is called the \emph{page} of the open
book. The \emph{monodromy} of an open book $(L,f)$ is given by the
return map of a flow transverse to the pages (all diffeomorphic to
$S$) and meridional near the binding, which is an element $h \in
Aut(S,\partial S)$, the group of (isotopy classes of)
diffeomorphisms of $S$ which restrict to the identity on $\partial
S$ . The group $Aut(S,\partial S)$ is also said to be the mapping
class group of $S$, and denoted by $\Gamma(S)$.

\medskip \noindent
An open book can also be described as follows. First consider the
mapping torus $$S(h)= [0,1]\times S/(1,x)\sim (0, h(x))$$ where $S$
is a compact oriented surface with $n=|\partial S|$ boundary
components and $h$ is an element of $Aut(S,\partial S)$ as above.
Since $h$ is the identity map on $\partial S$, the boundary
$\partial S(h)$ of the mapping torus $S(h)$ can be canonically
identified with $n$ copies of $T^2 = S^1 \times S^1$, where the
first $S^1$ factor is identified with $[0,1] / (0\sim 1)$ and the
second one comes from a component of $\partial S$. Now we glue in
$n$ copies of $D^2\times S^1$ to cap off $S(h)$ so that $\partial
D^2$ is identified with $S^1 = [0,1] / (0\sim 1)$ and the $S^1$
factor in $D^2 \times S^1$ is identified with a boundary component
of $\partial S$. Thus we get a closed $3$-manifold
$$M=M_{(S,h)}:= S(h) \cup_{n} D^2 \times S^1 $$
\emph{equipped with} an open book decomposition $(S,h)$ whose
binding is the union of the core circles in the $D^2 \times S^1$'s
that we glue to $S(h)$ to obtain $M$. See \cite{Gd}, \cite{Et2} for details.


\subsection{Legendrian knots and contact surgery }
\noindent A \emph{Legendrian knot} $K$ in a contact $3$-manifold $(M,\xi )$ is
a knot that is everywhere tangent to $\xi$. Any Legendrian knot
comes with a canonical {\em contact framing} (or
\emph{Thurston-Bennequin framing}), which is defined by a vector
field along $K$ that is transverse to $\xi$. We call $(M,\xi )$ (or
just $\xi$) {\em overtwisted} if it contains an embedded disc $D
\approx D^2\subset M$ with boundary $\partial D \approx S^1$ a
Legendrian knot whose contact framing equals the framing it receives
from the disc $D$. If no such disc exists, the contact structure
$\xi$ is called {\em tight}. Also if a contact $3$-manifold
$(M,\xi)$ is the boundary of a Stein manifold (resp. a symplectic
manifold) with certain compatibility conditions satisfied, then
$\xi$ is called \emph{Stein (holomorphically) fillable} (resp.
\emph{symplectically fillable}). See \cite{Et2} or \cite{OS} for the
complete definitions of fillability, and related facts. We will use the
following fact later.

\begin{thm} [\cite{EG}] \label{Eliashberg-Gromov}
Any symplectically fillable contact structure is tight.\\
$($ $\Rightarrow$ Any holomorphically fillable contact structure is
tight. $)$
\end{thm}

\noindent For any $p,q \in \mathbb{Z}$, a contact $(r)$-surgery
$(r=p/q)$ along a Legendrian knot $K$ in a contact manifold $(M,\xi
)$ was first described in \cite{DG}. It was proved in \cite{Ho} that
if $r=1/k$ with $k\in \mathbb{Z}$, then the resulting contact
structure is unique up to isotopy. In particular, a contact
$\pm1$-surgery along a Legendrian knot $K$ on a contact manifold
$(M,\xi)$ determines a unique surgered contact manifold which will
be denoted by $(M,\xi)_{(K,\pm1)}$.

\medskip \noindent The most general result along these lines is:
\begin{thm} [\cite{DG}] \label{DingGeigesmain}
Every closed contact $3$-manifold $(M, \xi )$ can be
obtained via contact $(\pm 1)$-surgery on a Legendrian link in
$(S^3, \xi _{st})$.
\end{thm}

\medskip \noindent
Any closed contact $3$-manifold $(M,\xi)$ can be described by a
\emph{contact surgery diagram} drawn in $(\mathbb{R}^3, \xi_0) \subset
(S^3,\xi_{st})$. By Theorem \ref{DingGeigesmain}, there is
a contact surgery diagram for $(M,\xi)$ such that the contact
surgery coefficient of any Legendrian knot in the diagram is $\pm1$.
For any oriented Legendrian knot $K$ in $(\mathbb{R}^3, \xi_0)$, we
compute the \emph{Thurston-Bennequin number} $tb(K)$, and the
\emph{rotation number} $rot(K)$ as
\begin{center}
$tb(K)=bb(K) - (\#$ of $\mbox{left cusps of K})$, \\
$rot(K)=\displaystyle{\frac{1}{2}} [(\#$ of downward cusps$)- (\#$
of upward cusps$)]$
\end{center}
where $bb(K)$ is the blackboard framing of $K$.

\medskip \noindent If a contact surgery diagram for $(M,\xi)$ is given,
we can also get the smooth surgery diagram for the underlying $3$-manifold $M$.
Indeed, for a Legendrian knot $K$ in a contact surgery diagram, we have:
\medskip
\begin{center}
Smooth surgery coefficient of $K$ = Contact surgery coefficient of
$K$ $+ \; tb(K)$
\end{center}
\medskip
\noindent For more details see \cite{OS} and \cite{Gm}.


\subsection{Compatibility and stabilization}

\noindent A contact structure $\xi$ on a $3$-manifold $M$ is said to be
\emph{supported by an open book} $(L,f)$ if $\xi$ is isotopic
to a contact structure given by a $1$-form $\alpha$ such that
\begin{enumerate}
\item $d\alpha$ is a positive area form on each page $S\approx f^{-1}($pt$)$ of the open book and
\item $\alpha>0$ on $L$ (Recall that $L$ and the pages are oriented.)
\end{enumerate}

\noindent When this holds, we also say that the open book $(L,f)$ is
\emph{compatible with the contact structure} $\xi$ on $M$.

\begin{defn} \label{stabilization}
A positive (resp., negative) stabilization $S^{+}_{K}(S,h)$ (resp.,
$S^{-}_K(S,h)$) of an abstract open book $(S,h)$ is the open book
\begin{enumerate}
\item with page $S'=S \cup \text{ 1-handle}$ and
\item monodromy $h'=h \circ D_K$ (resp., $h'=h \circ D_K^{-1}$)
where $D_K$ is a right-handed Dehn twist along a curve $K$ in $S'$
that intersects the co-core of the 1-handle exactly once.
\end{enumerate}
\end{defn}

\medskip \noindent
Based on the result of Thurston and Winkelnkemper \cite{TW} which
introduced open books into the contact geometry, Giroux proved the
following theorem strengthening the link between open books and
contact structures.

\begin{thm} [\cite{Gi}] \label{Giroux}
Let $M$ be a closed oriented $3$-manifold. Then there is a
one-to-one correspondence between oriented contact structures on $M$
up to isotopy and open book decompositions of $M$ up to positive
stabilizations: Two contact structures supported by the same open
book are isotopic, and two open books supporting the same contact
structure have a common positive stabilization.
\end{thm}

\noindent Following fact was first implied in \cite{LP}, and then in
\cite{AO}. The given version below is due to Giroux and Matveyev.
For a proof, see \cite{OS}.

\begin{thm}\label{Holofill=PosDehnTwist}
A contact structure $\xi$ on $M$ is holomorphically fillable if and
only if  $\xi$ is supported by some open book whose monodromy admits
a factorization into positive Dehn twists only.
\end{thm}

\noindent For a given fixed open book $(S,h)$ of a $3$-manifold $M$,
there exists a unique compatible contact structure up to isotopy on
$M=M_{(S,h)}$ by Theorem \ref{Giroux}. We will denote this contact
structure by $\xi_{(S,h)}$. Therefore, an open book $(S,h)$
determines a unique contact manifold $(M_{(S,h)},\xi_{(S,h)})$ up to
contactomorphism.

\medskip
\noindent Taking a positive stabilization of $(S,h)$ is actually
taking a special Murasugi sum of $(S,h)$ with the positive Hopf
band $(H^+,D_{\gamma})$ where $\gamma \subset H^+$ is the core circle. Taking a Murasugi sum of two open books corresponds to taking
the connect sum of $3$-manifolds associated to the open books. The
proofs of the following facts can be found in \cite{Gd}, \cite{Et2}.
\begin{thm} \label{connectsummingwithS^3}
$(M_{S^{+}_{K}(S,h)},\xi_{S^{+}_{K}(S,h)}) \cong
(M_{(S,h)},\xi_{(S,h)}) \# (S^3,\xi_{st}) \cong
(M_{(S,h)},\xi_{(S,h)}).$
\end{thm}


\begin{thm} \label{surgeryonbook}
Let $(S,h)$ be an open book supporting the contact manifold
$(M,\xi).$ If $K$ is a Legendrian knot on the page $S$ of the open
book, then
\[
(M,\xi)_{(K,\pm 1)}= (M_{(S,h\circ D_K^{\mp})}, \xi_{(S,h \circ
D_K^{\mp})}).
\]
\end{thm}


\subsection{Homotopy invariants of contact structures}
\noindent The set of oriented $2-$plane fields on a given $3$-manifold $M$ is
identified with the space $Vect(M)$ of nonzero vector fields on $M$.
$v_1, v_2 \in Vect(M)$ are called \emph{homologous} (denoted by
$v_1\sim v_2$) if $v_1$ is homotopic to $v_2$ in $M \setminus B$ for
some $3-$ball $B$ in $M$. The space $Spin^c(M)$ of all spin$^c$
structures on $M$ is the defined to be the quotient space $Vect(M)/
\sim$. Therefore, any contact structure $\xi$ on $M$ defines a
spin$^c$ structure $\t_{\xi} \in Spin^c(M)$ which depends only on
the homotopy class of $\xi$. As the first invariant of $\xi$, we will use
the first Chern class $c_1(\xi) \in H^2(M;\Z)$
(considering $\xi$ as a complex line bundle on $M$). For a
spin$^c$ structure $\t_{\xi}$, whose first Chern class
$c_1(\t_{\xi})(:= c_1(\xi))$ is torsion, the obstruction to homotopy
of two $2$-plane fields (contact structures) both inducing
$\t_{\xi}$ can be captured by a single number. This obstruction is
the \emph{$3$-dimensional invariant} $d_3(\xi)$ of $\xi $). To
compute $d_3(\xi)$, suppose that a compact almost complex 4-manifold
$(X,J)$ is given such that $\partial X=M$, and $\xi $ is the complex
tangencies in $TM$, i.e., $\xi =TM\cap J(TM)$. Let $\sigma (X),\chi
(X)$ denote the signature and Euler characteristic of $X$,
respectively. Then we have

\begin{thm}[\cite{Gm}]\label{d3_original}
If $c_1(\xi)$ is a torsion class, then the rational number
\[ d_3(\xi)=\frac{1}{4}\bigl( c_1^2(X,J)-3\sigma (X)-2\chi (X)\bigr)\]
is an invariant of the homotopy type of the $2$-plane field $\xi$.
Moreover, two $2$-plane fields $\xi _1$ and $\xi _2$ with
$\mathfrak{t}_{\xi_1}=\t_{\xi_2}$ and  $c_1(\t _{\xi
_i})=c_1(\xi_i)$ a torsion class are homotopic if and only if
$d_3(\xi _1)=d_3(\xi _2).$ \qed
\end{thm}

\noindent As a result of this fact, if $(M,\xi)$ is given by a
contact $\pm1$-surgery on a link, then we have

\begin{cor}
[\cite{DGS}] \label{d3_computation} Suppose that $(M, \xi )$, with
$c_1(\xi )$ torsion, is given by a contact $(\pm 1)$-surgery on a
Legendrian link ${\mathbb {L}}\subset (S^3, \xi _{st})$ with
$tb(K)\neq 0$ for each $K\subset {\mathbb L}$ on which we perform
contact $(+1)$-surgery. Let $X$ be a $4$-manifold such that
$\partial X=M$. Then
$$d_3(\xi )=\frac{1}{4}\bigl( c^2-3\sigma (X)-2\chi (X)\bigr) +s,$$
where $s$ denotes the number of components in ${\mathbb {L}}$ on
which we perform $(+1)$-surgery, and $c\in H^2 (X; \Z )$ is the
cohomology class determined by $c(\Sigma_K)=rot(K)$ for each
$K\subset \mathbb {L}$, and $[\Sigma_K]$ is the homology class in
$H_2(X)$ obtained by gluing the Seifert surface of $K$ with the core
disc of the $2$-handle corresponding $K$.
\end{cor}

\noindent We use the above formula as follows: Suppose $\mathbb{L}$
has $k$ components. Write $\mathbb{L}={\sqcup_i}^k K_i$. By
converting all contact surgery coefficients to the topological ones,
and smoothing each cusp in the diagram, we get a framed link (call
it $\mathbb{L}$ again) describing a simply connected $4$-manifold
$X$ such that $\partial X= M$. Using this description, we compute
\begin{center}
$\chi(X)=1+ k$, and $\sigma(X)=\sigma(\mathcal{A}_{\mathbb{L}})$
\end{center}
where $\mathcal{A}_{\mathbb{L}}$ is the linking matrix of
$\mathbb{L}$. Using the duality, the number $c^2$ is computed as
\begin{center}
$c^2=(PD (c))^2=[b_1 \; b_2 \cdots b_k]\mathcal{A}_{\mathbb{L}}[b_1
\; b_2 \cdots b_k]^T$
\end{center}
where $PD(c) \in H_2(X,\partial X; \Z)$ is the Poincar\'e dual of $c$,
the row matrix $[b_1 \; b_2 \cdots b_k]$ is the unique solution to
the linear system
\begin{center}
$\mathcal{A}_{\mathbb{L}}[b_1 \; b_2 \cdots b_k]^T=[rot(K_1) \;
rot(K_2) \cdots rot(K_k)]^T$.
\end{center}
Here the superscript ``$\,T\,$" denotes the transpose operation in
the space of matrices. See \cite{DGS}, \cite{Gm} for more details.


\subsection{Right-veering diffeomorphisms}
\begin{figure}[ht]
  \begin{center}
   \includegraphics{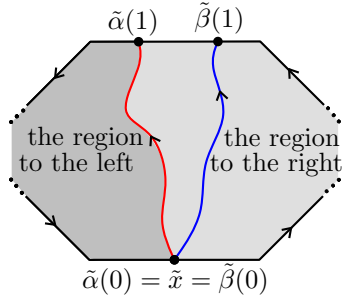}
   \caption{Lifts of $\alpha$ and $\beta$ in the universal cover $\tilde S$.}
  \label{BetaRightAlfa}
    \end{center}
\end{figure}

\noindent For a given compact oriented surface $S$ with nonempty boundary
$\partial S$, let $Dehn^{+}(S,\partial S) \subset Aut(S,\partial S)$
be the submonoid of product of all positive Dehn twists. In
\cite{HKM}, another submonoid $Veer(S,\partial S)$ of all \rv
elements in $Aut(S,\partial S)$ was introduced and studied. They
defined \rv elements of $Aut(S,\partial S)$ as follows: Let $\alpha$
and $\beta$ be isotopy classes (relative to the endpoints) of
properly embedded oriented arcs $[0,1]\rightarrow S$ with a common
initial point $\alpha(0)=\beta(0)=x\in \partial S$. Let $\pi:\tilde
S\rightarrow S$ be the universal cover of $S$ (the interior of
$\tilde S$ will always be $\R^2$ since $S$ has at least one boundary
component), and let $\tilde x\in \partial \tilde S$ be a lift of
$x\in \partial S$. Take lifts $\tilde \alpha$ and $\tilde \beta$ of
$\alpha$ and $\beta$ with $\tilde \alpha(0) =\tilde \beta(0)=\tilde
x$. $\tilde \alpha$ divides $\tilde S$ into two regions -- the
region ``to the left'' (where the boundary orientation induced from
the region coincides with the orientation on $\tilde \alpha$) and
the region ``to the right''. We say that $\beta$ is {\em to the
right} of $\alpha$ if either $\alpha=\beta$ (and hence
$\tilde\alpha(1)=\tilde\beta(1)$), or $\tilde\beta(1)$ is in the
region to the right (see Figure \ref{BetaRightAlfa}).

\noindent Alternatively, isotop $\alpha$ and $\beta$, while fixing
their endpoints, so that they intersect transversely (this include
the endpoints) and with the fewest possible number of intersections.
Assume that $\alpha\not=\beta$.  Then in the universal cover $\tilde
S$, $\tilde\alpha$ and $\tilde\beta$ will meet only at $\tilde x$.
If not, subarcs of $\tilde \alpha$ and $\tilde\beta$ would cobound a
disk $D$ in $\tilde S$, and we could use an innermost disk argument
on $\pi(D)\subset S$ to reduce the number of intersections of
$\alpha$ and $\beta$ by isotopy. Then $\beta$ is to the right of
$\alpha$ if $int(\tilde\beta)$ lies in the region to the right. As
an alternative to passing to the universal cover, we simply check to
see if the tangent vectors $(\dot \beta(0), \dot\alpha(0))$ define
the orientation on $S$ at $x$.

\begin{defn}
Let $h:S \to S$ be a diffeomorphism that restricts to the identity
map on $\partial S$. Let $\alpha$ be a properly embedded oriented
arc starting at a basepoint $x \in \partial S$. Then $h$ is {\em
right-veering} (that is, $h \in Veer(S,\partial S)$) if for every
choice of basepoint $x \in
\partial S$ and every choice of $\alpha$ based at $x$, $h(\alpha)$
is to the right of $\alpha$ (at $x$). If $C$ is a boundary component
of $S$, we say is $h$ is {\em right-veering with respect to $C$} if
$h(\alpha)$ is to the right of $\alpha$ for all $\alpha$ starting at
a point on $C$.
\end{defn}
\noindent It turns out that $Veer(S,\partial S)$ is a submonoid and
we have the inclusions:
\medskip
\begin{center}
 $Dehn^{+}(S,\partial S) \subset Veer(S,\partial S) \subset Aut(S,\partial S)$.
\end{center}
\medskip
\noindent In \cite{HKM}, they proved the following theorem which is
hard to use but still can be used to distinguish tight structure in
some cases.
\begin{thm}[\cite{HKM}] \label{hokama}
A contact structure $(M,\xi)$ is tight if and only if all of its
compatible open book decompositions $(S,h)$ have right-veering $h
\in Veer(S,\partial S) \subset Aut(S,\partial S)$.
\end{thm}

\section{The proofs of results} \label{Results}

\noindent We first prove that the submonoids \D and \V are actually
the same in our particular case.
\begin{lem}  \label{Dehn=Veer}
\D = \V for the surface $\Sigma$ given in Figure \ref{our surface}.
\end{lem}
\begin{proof}
The inclusion $Dehn^{+}(S,\partial S) \subset Veer(S,\partial S)$ is
true for a general compact oriented surface $S$ with boundary (see
Lemma 2.5. in \cite{HKM} for the proof). Now, suppose that $\phi
\in$ \V $\subset$ \au. Then we can write $\phi$ in the form
\begin{center}
$\phi={D_a}^p{D_b}^q{D_c}^r$ for some $p,q,r \in \Z$.
\end{center}
We will show that $p\geq0, q\geq0, r\geq0$. Consider the properly
embedded arc $\alpha \subset \Sigma$ one of whose end points is
$x \in \partial \Sigma$ as shown in the Figure \ref{dehnveer}. Note
that, for any $p, q, r \in \Z$, ${D_c}^r$ fixes $\alpha$, and also
any image ${D_a}^p{D_b}^q(\alpha)$ of $\alpha$ because $c$ does not
intersect any of these arcs.
\begin{figure}[ht]
  \begin{center}
   \includegraphics{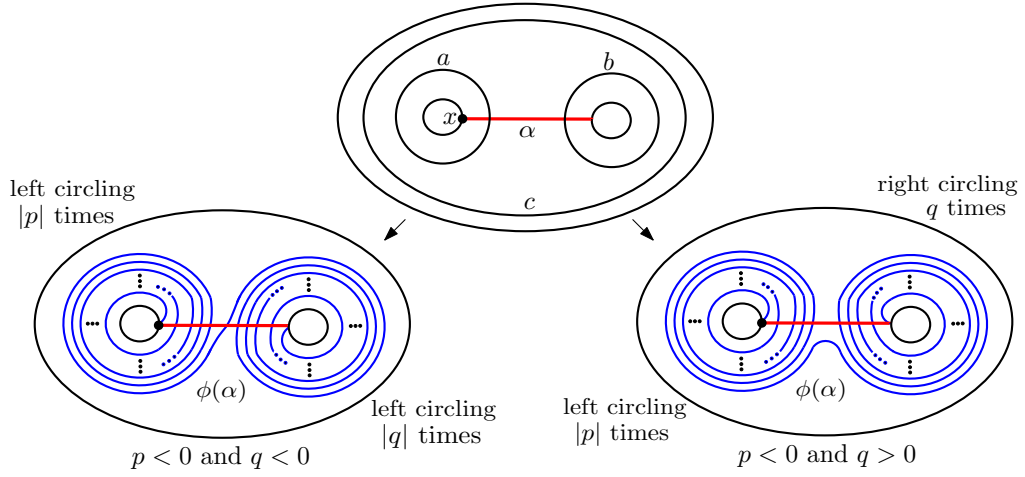}
   \caption{The arc $\alpha$ and its image $\phi(\alpha)={D_a}^p{D_b}^q(\alpha)$. }
  \label{dehnveer}
    \end{center}
\end{figure}
Assume at least one of $p, q,$ or $r$ is strictly negative. First
assume that $p<0$. Then consider two possible different images
$\phi(\alpha)={D_a}^p{D_b}^q(\alpha)$ of $\alpha$ corresponding to
whether $q<0$ or $q>0$ (See Figure \ref{dehnveer}). Since we are not
allowed to rotate any boundary component, clearly $\phi(\alpha)$ is
to the left of $\alpha$ at the boundary point $x$. Equivalently,
$\phi(\alpha)$ is not to the right of $\alpha$ at $x$ which implies
that $h$ is not right-veering with respect to the boundary component
parallel to $a$. Therefore, $\phi \notin$ \V which is a
contradiction. Now by symmetry, we are also done for the case $q<0$.
Finally, exactly the same argument (with a different choice of arc
one of whose end points is on the boundary component parallel to
the curve $c$) will work for the case when $r<0$.
\end{proof}

\begin{lem}  \label{Tight=holofil}
Let $(M,\xi)$ be a contact manifold. Assume that $\xi$ is supported
by $(\Sigma,\phi)$ where $\phi \in$ \au. Then $\xi$ is tight if and
only if $\xi$ is holomorpfically fillable.
\end{lem}

\begin{proof}
Assume that $\xi$ is tight. Since $\phi \in$ \au, there exists
integers $p, q, r$ such that $\phi={D_a}^p{D_b}^q{D_c}^r$. As $\xi$
is tight, the monodromy of any open book supporting $\xi$ is \rv by
Theorem \ref{hokama}. In particular, we have $\phi \in$ \V since
$(\Sigma, \phi)$ supports $\xi$. Therefore, $\phi \in$ \D by Lemma \ref{Dehn=Veer},
and so $p\geq0, q\geq0, r\geq0$. Thus, $\xi$ is
holomorphically fillable by Theorem \ref{Holofill=PosDehnTwist}.
Converse statement is a consequence of Theorem
\ref{Eliashberg-Gromov}.
\end{proof}

\noindent Now, the following corollary of Lemma
\ref{Tight=holofil} is immediate:
\begin{cor}  \label{Tight=positiveDenhtwists}
Let $(M,\xi)$ be a contact manifold. Assume that $\xi$ is supported
by $(\Sigma,\phi)$ where $\phi \in$ \au. Then
\begin{center}
 $\xi$ is tight $\Longleftrightarrow$ $\phi={D_a}^p{D_b}^q{D_c}^r$ with $p\geq0, q\geq0,
 r\geq0$.
\end{center}
\end{cor}
\vspace{-.75cm} \hspace{12.7cm} $\square$

\medskip
\begin{proof}[\textbf{Proof of Theorem \ref{Contact_Type}}]
Let $(M,\xi)$ be a contact manifold supported by the open book
$(\Sigma,\phi_{p,q,r})$ where $\phi_{p,q,r}={D_a}^p{D_b}^q{D_c}^r
\in$ \A for $p,q,r \in \Z$. As explained in \cite{EO},
$(M,\xi)=(M_{(\Sigma,\phi_{p,q,r})},\xi_{(\Sigma,\phi_{p,q,r})})$ is
given by the contact surgery diagram in Figure
\ref{contactstructure}. Then we apply the algorithm given in
\cite{DG} and \cite{DGS} to convert each rational coefficient into
$\pm1$'s, and obtain the diagram given in Figure
\ref{modelstructure}.

\begin{figure}[ht]
  \begin{center}
   \includegraphics{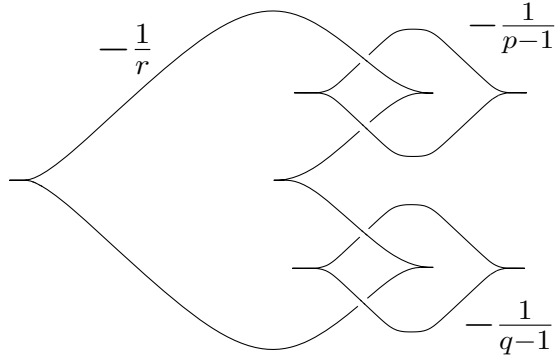}
   \caption{Contact surgery diagram corresponding to $(\Sigma,\phi_{p,q,r})$.}
  \label{contactstructure}
    \end{center}
\end{figure}

\begin{figure}[ht]
  \begin{center}
   \includegraphics{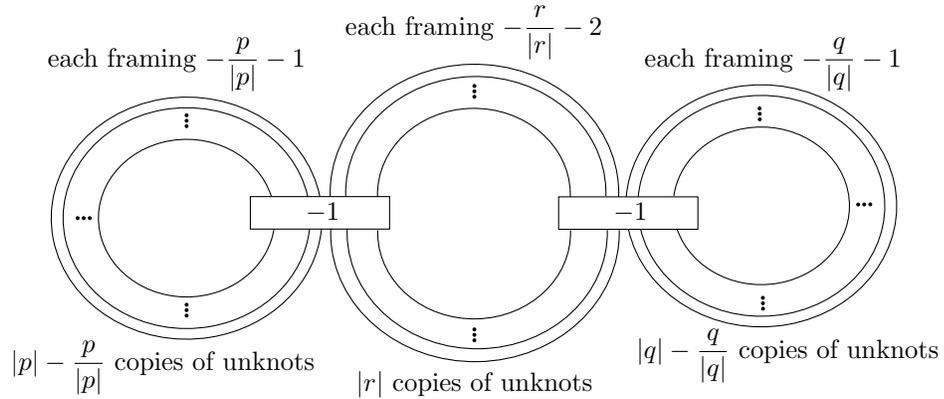}
   \caption{Smooth surgery diagram corresponding to Figure \ref{modelstructure}. }
  \label{smoothmodel}
    \end{center}
\end{figure}

\medskip
\noindent To determine the topological (or smooth) type of
$(M,\xi)$, we start with the diagram in Figure \ref{modelstructure}.
Then by converting the contact surgery coefficients into the smooth
surgery coefficients, we get the corresponding smooth surgery
diagram in Figure \ref{smoothmodel} where each curve is an unknot.

\begin{figure}[ht]
  \begin{center}
   \includegraphics{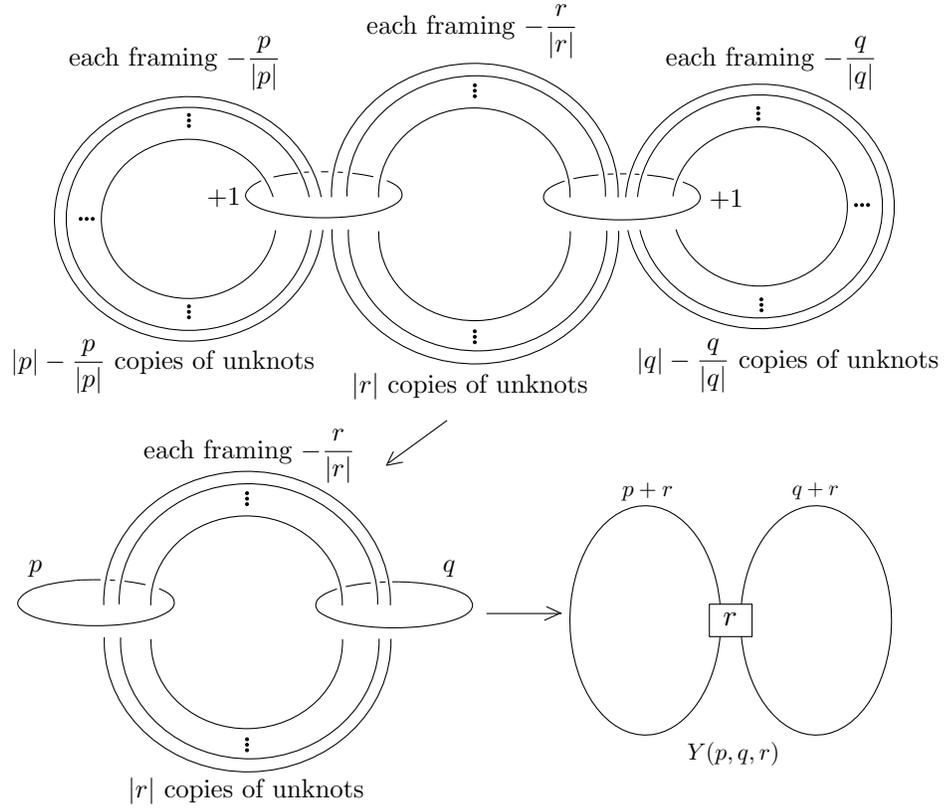}
   \caption{Squence of blow-ups and blow-downs.}
  \label{smoothmodel2}
    \end{center}
\end{figure}

\medskip \noindent Now we modify this diagram using a sequence of blow-ups and
blow-downs. These operations do not change smooth type of $M$. We
first blow up the diagram twice so that we unlink two $-1$ twists.
Then we blow down each unknot in the most left and the most right
families. Finally we blow down each unknot of the family in the
middle. We illustrate these operations in Figure \ref{smoothmodel2}.
To keep track the surgery framings, we note that each blow-up
increases the framing of any unknot by 1 if the unknot passes
through the corresponding twist box in Figure \ref{smoothmodel}. So
we get the first diagram in Figure \ref{smoothmodel2}. Blowing each
member down on the left (resp. right) decreases the framing of the
left (resp. right) $+1$-unknot by $-\frac{p}{|p|}$ (resp.
$-\frac{q}{|q|}$).  Since there are $|p|-\frac{p}{|p|}$ blow-downs
on the left and $|q|-\frac{q}{|q|}$ blow-downs on the right, we get
the second diagram in Figure \ref{smoothmodel2}. Finally, if we blow
down each ($-\frac{r}{|r|}$)-unknot in the middle family, we get the
last diagram. Note that each blow-down decreases the framing by
$-\frac{r}{|r|}$, and introduces a $\frac{r}{|r|}$ full twist.
Hence, we showed that $(M,\xi)$ is contactomorphic to
$(Y(p,q,r),\xi_{p,q,r})$. The statements (1) and (2) are the
consequences of Corollary \ref{Tight=positiveDenhtwists}.
\end{proof}

\noindent We now examine the special case where $Y(p,q,r)$ is
homeomorphic to $3$-sphere $S^3$. The following lemma lists all
planar contact structures on $S^3$ with binding number less than or
equal to three.

\begin{lem} \label{Case_of_Sphere} Suppose that $(Y(p,q,r), \xi_{p,q,r})$
is contactomorphic to $(S^3, \xi)$ for some contact structure $\xi$ on $S^3$.
Then Table \ref{Structures_on_Sphere} lists all possible values of $(p, q, r)$,
the corresponding $\xi$ (in terms of the $d_3$-invariant), and its binding number.
\end{lem}

\begin{table}[ht]
  \begin{center}
   \includegraphics{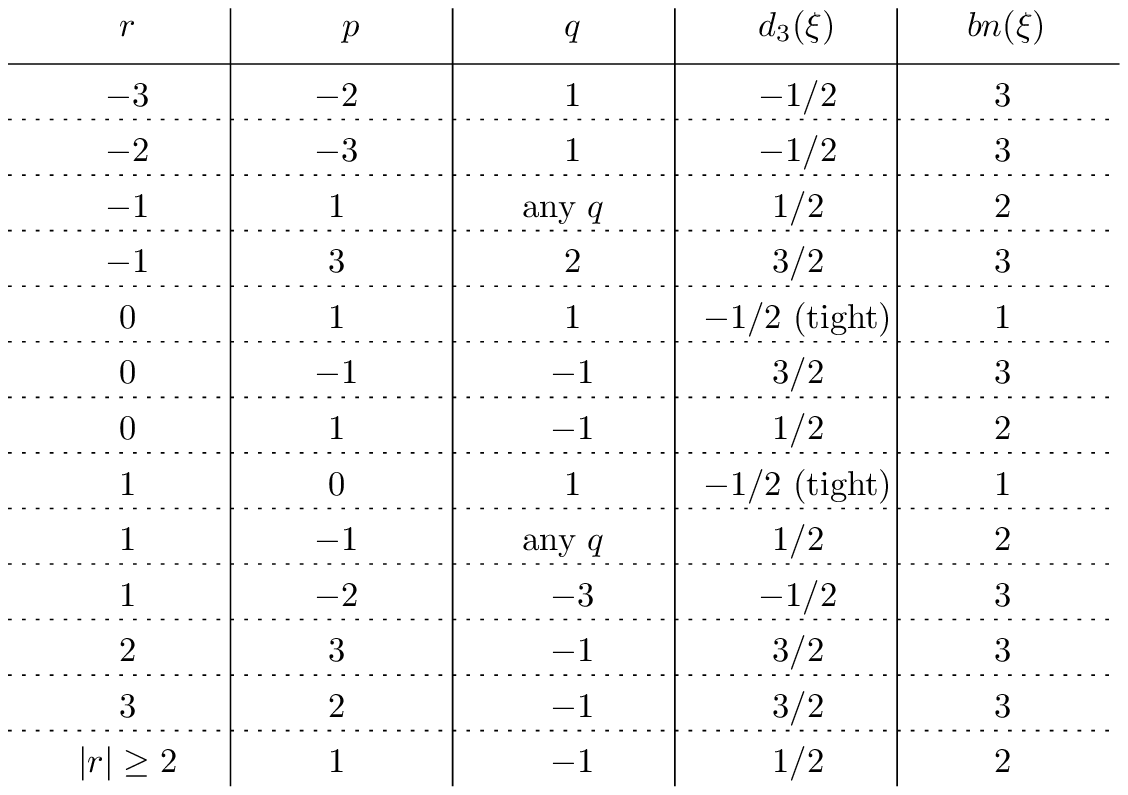}
      \caption{All planar contact structures on $S^3$ with binding number $\leq 3$.}
  \label{Structures_on_Sphere}
    \end{center}
\end{table}

\begin{proof} The proof is the direct consequence of the discussion given in the
proof of Lemma 5.5 in \cite{EO}. We remark that the interchanging
$p$ and $q$ does not affect the contact structure in Figure
\ref{modelstructure}, so we do not list the possibilities for
$(p,q,r)$ that differ by switching $p$ and $q$. Note that in Table \ref{Structures_on_Sphere}
there are only two contact structures (up to isotopy) on
$S^3$ with binding number 3, namely, the ones with $d_3$-invariants
$-1/2$ and $3/2$.
\end{proof}

\begin{proof}[\textbf{Proof of Theorem \ref{the_List}}]
We will use the results of Theorem \ref{bindingnumber=<2}, Theorem
\ref{Contact_Type}, and Lemma \ref{Case_of_Sphere}. Consider the
$3$-sphere $S^3$ in Theorem \ref{bindingnumber=<2} as the lens space
$L(1,\pm1)$. By Theorem \ref{bindingnumber=<2}, for any contact
manifold $(Y,\eta)$ with $sg(\eta)=0$ and $bn(\eta)\leq2$, we have
either

\begin{enumerate}
\item $(Y,\eta)\cong (S^3, \xi_{st})$ if $bn(\eta)=1$,
\item $(Y,\eta)\cong (L(m,-1), \eta_m)$ for some $m \geq 2$ if $bn(\eta)=2$, and $\eta$ is tight,
\item $(Y,\eta)\cong (L(m,1), \eta_m)$ for some $m\geq0$ if $bn(\eta)=2$, and $\eta$ is overtwisted (for $m\neq0$).
\end{enumerate}

\noindent where $\eta_m$ is the contact structure on the lens space
$L(m,-1)$ (or $L(m,1)$) given by the contact surgery diagram
consisting of a single family of Legendrian unknots (with
Thurston-Benequen number $-1$) such that each member links all the
other members of the family once, and each contact surgery
coefficient is $-1$ (if $\eta_m$ is tight) or $+1$ (if $\eta_m$ is
overtwisted). These are illustrated by the diagrams $(\ast)$ and
$(\star)$ in Figure \ref{lensspaces}, respectively. Notice the
exceptional cases: $m=1$ in $(\ast)$, and $m=0$ in $(\star)$.

\begin{figure}[ht]
  \begin{center}
   \includegraphics{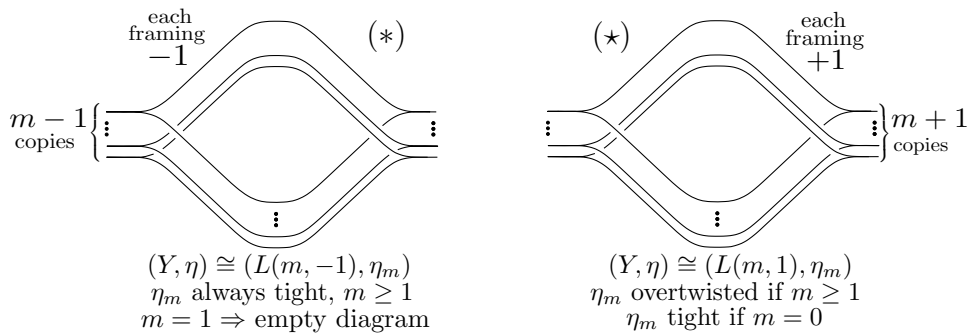}
   \caption{Contact surgery diagrams for $(Y,\eta).$}
  \label{lensspaces}
    \end{center}
\end{figure}

\medskip \noindent Now, if $(M,\xi)$ is a contact manifold with $sg(\xi)=0$
and $bn(\xi)=3$, then by the definitions of these invariants there
exists an open book $(\Sigma,\phi)$ supporting $\xi$. Therefore, by
Theorem \ref{Contact_Type}, $(M,\xi)$ is contactomorphic to
$(Y(p,q,r),\xi_{p,q,r})$ for some $p,q,r \in \Z$, and the contact
surgery diagram of $\xi$ is given in Figure \ref{modelstructure}.
However, $p,q,r$ can not be arbitrary integers because there are
several cases where the diagram in Figure \ref{modelstructure}
reduces to either $(\ast)$ or
$(\star)$ in Figure \ref{lensspaces} for some $m$. So for those
values of $p,q,r$, $(M,\xi)$ can not be contactomorphic to
$(Y(p,q,r),\xi_{p,q,r}) \cong (Y,\eta)$ because $bn(\xi)=3\neq2\geq
bn(\eta)$. Therefore, we have to determine those cases.

\medskip \noindent If $|p\,|\geq2$ and $|q|\geq2$, then the only triples
$(p,q,r)$ giving $L(m,\pm1)'s$ are $(-2,q,1)$ and $(2,q,-1)$. Furthermore,
if we assume also that $|r|>1$, then the Seifert fibered manifolds $Y(p,q,r)$
are not homeomorphic to even a lens space $L(m,n)$ for any $m,n$ (for instance,
see Chapter 5 in \cite{Or}). As a result, we immediately obtain $bn(\xi_{p,q,r})=3$
for $|p\,|\geq2$ and $|q|\geq2$ and $|r|\geq2$. Therefore, to finish the proof of
the theorem, it is enough to analyze the cases where $|p\,|<2$ or $|q|<2$, and
the cases $(-2,q,1)$ and $(2,q,-1)$ for any $q$. As we remarked before, we do
not need to list the possibilities for $(p,q,r)$ that differ by switching $p$
and $q$. We first consider $r=0$, $\pm1$, $\pm2$, and then the cases $r>2$ and
$r<-2$. In Table \ref{r_is_0} - \ref{r_is_less_-2}, we list all possible
$(M,\xi)$ for each of these cases.

\begin{rem} \label{Main_Remark}
To determine the binding number $bn(\xi)$ in any row of any table
below, we simply first check the topological type of the manifold
under consideration. If $M \approx S^3$, we determine the
corresponding binding number using Table \ref{Structures_on_Sphere}.
If the topological type is not $L(m,1)$ or $L(m,-1)$, then we
immediately get that $bn(\xi)=3$. If $M\approx L(m,1)$ with $m>1$,
then we first compute $c_1(\xi)$. If $c_1(\xi)\neq0$, then
$bn(\xi)=3$ as $c_1(\eta_m)=0$ for any $\eta_m$ given above. If
$c_1(\xi)$=0, we compute the $d_3(\xi)$ using the $4$-manifold
defined by the surgery diagram in Figure \ref{smoothmodel}. (Indeed,
we can use the formula for $d_3$ given in Corollary
\ref{d3_computation} as long as $c_1(\xi)$ is torsion. In
particular, whenever $H^2(M)$ is finite, then $d_3$ is computable).
Then if $d_3(\xi)=d_3(\eta_m)=(-m+3)/4$, then $\xi$ is isotopic to
$\eta_m$ which implies that $bn(\xi)= 2$ by Theorem
\ref{bindingnumber=<2}. Otherwise $bn(\xi)=3$. In the case that
$M\approx L(m,-1)$ with $m>1$, we first ask if $\xi$ is tight. If it
is tight (which is the case if and only if $p \geq 0, q \geq 0, r \geq 0$),
then $bn(\xi)=2$ (again by Theorem \ref{bindingnumber=<2}) since the
tight structure on $L(m,-1)$ is unique (upto isotopy). If it is
overtwisted (which is the case if and only if at least one of $p, q, r$ is
negative), then $bn(\xi)=3$ because $\xi$ is not covered in Theorem
\ref{bindingnumber=<2}. As a final remark, sometimes the contact
structure $\xi$ can be viewed as a positive stabilization of some
$\eta_m$. For these cases we immediately obtain that $bn(\xi)=2$
because positive stabilizations do not change the isotopy classes of
contact structures.
\end{rem}

\noindent To compute the $d_3$-invariant of $\xi_{p,q,r}$
(for $c_1(\xi_{p,q,r})$ torsion), we will use the $(n+1)\times(n+1)$
matrices $A_n$ ($n\geq1$), $B_n$ ($n\geq1$), and $C_n$ ($n\geq4$) given below.
It is a standard exercise to check that
\begin{enumerate}
\item $\sigma(A_n)=n-1$ if $n\geq1$, and $\sigma(C_n)=n-1$ if $n\geq4$. \\
\item $\sigma(B_n)=n-3$ if $n\geq3$, and $\sigma(B_n)=0$ if $n=1,2$. \\
\item The system $A_n[\mathbf{b}]^T_{n+1}=[\mathbf{0}]^T_{n+1}$ has trivial solution $[\mathbf{b}]^T_{n+1}=[\mathbf{0}]^T_{n+1}$ where \\ $[\mathbf{b}]_{n+1}=[b_1 \; b_2 \cdots b_{n+1}]$, $[\mathbf{0}]_{n+1}=[0 \; 0 \cdots 0]$ are $(n+1) \times 1$ row matrices.
\end{enumerate}

\begin{table}[ht]
   \includegraphics{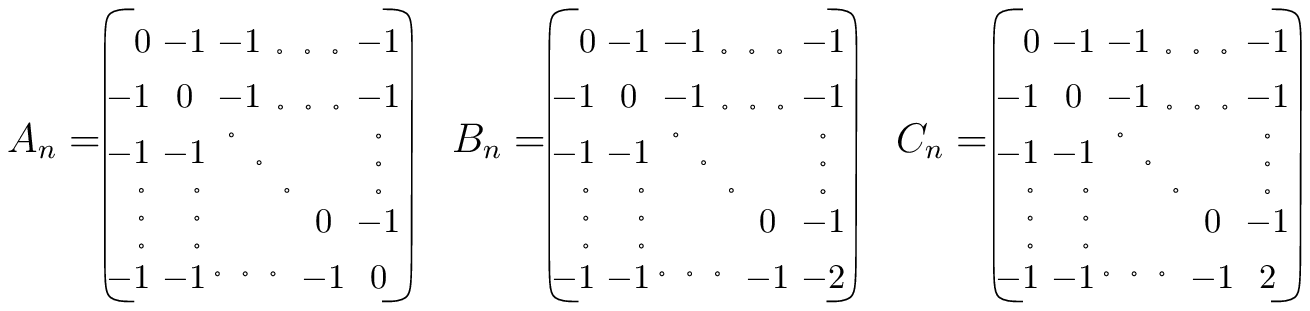}
\end{table}

\noindent In some cases, $A_n$ appears (as a block matrix) in the
linking matrix $\mathcal{L}_{p,q,r}$ of the framed link
$\mathbb{L}_{p,q,r}$ given in Figure \ref{smoothmodel}. On the other
hand, $B_n$ and $C_n$ are very handy  when we diagonalize
$\mathcal{L}_{p,q,r}$ to find its signature. As we discussed before,
the link $\mathbb{L}_{p,q,r}$ defines a $4$-manifold $X_{p,q,r}$
with $\partial X=M$. So we have

$$\begin{array}{crl}
\sigma(X_{p,q,r})&=&\sigma(\mathcal{L}_{p,q,r}), \\
\chi(X_{p,q,r})&=&1+ (\# \; \textrm{of components of} \; \mathbb{L}_{p,q,r}), \\
c^2&=&[\mathbf{b}]_k \mathcal{L}_{p,q,r}[\mathbf{b}]^T_k
\end{array}$$

\noindent where $[\mathbf{b}]^T_k$ is the solution to the linear
system $\mathcal{L}_{p,q,r}[\mathbf{b}]^T_k=[rot(K_1) \;rot(K_2)
\cdots rot(K_k)]^T$ with $K_1, K_2, \cdots K_k$ being the components
of $\mathbb{L}_{p,q,r}$.

\medskip \noindent To compute the first Chern class $c_1(\xi_{p,q,r})
\in H^2(M)$, note that in Figure \ref{modelstructure}, the rotation
number of any member in the family corresponding to $r$ is $\pm1$
(depending on how we orient them). We will always orient them so
that their rotation numbers are all $+1$. On the other hand, the
rotation number is $0$ for any member in the family corresponding to
$p$ and $q$. Therefore, $c_1(\xi_{p,q,r})=PD^{-1}(\mu_1+\mu_2 +
\cdots + \mu_{|r|})$ where $\mu_i \in H_1(M)$ is the class of the
meridian of the Legendrian knot $K_i$ in the family corresponding to
$r$. Then we compute $H_1(M)$ (which is isomorphic to  $H^2(M)$ by
Poincar\'e duality) as
\begin{center}
$H_1(M)=\langle \; \mu_1, \mu_2, \cdots, \mu_k | \;
\mathcal{L}_{p,q,r}[\mathbf{\mu}]^T_k=[\mathbf{0}]^T_k \; \rangle$
\end{center}
where $[\mu]_k=[\mu_1 \; \mu_2 \cdots \mu_k]$ is the $k \times 1$
row matrix. The final step is to understand
$PD(c_1(\xi_{p,q,r}))=\mu_1+\mu_2 + \cdots + \mu_{|r|}$ in this
presentation of $H_1(M)$.

\begin{table}[ht]
   \includegraphics{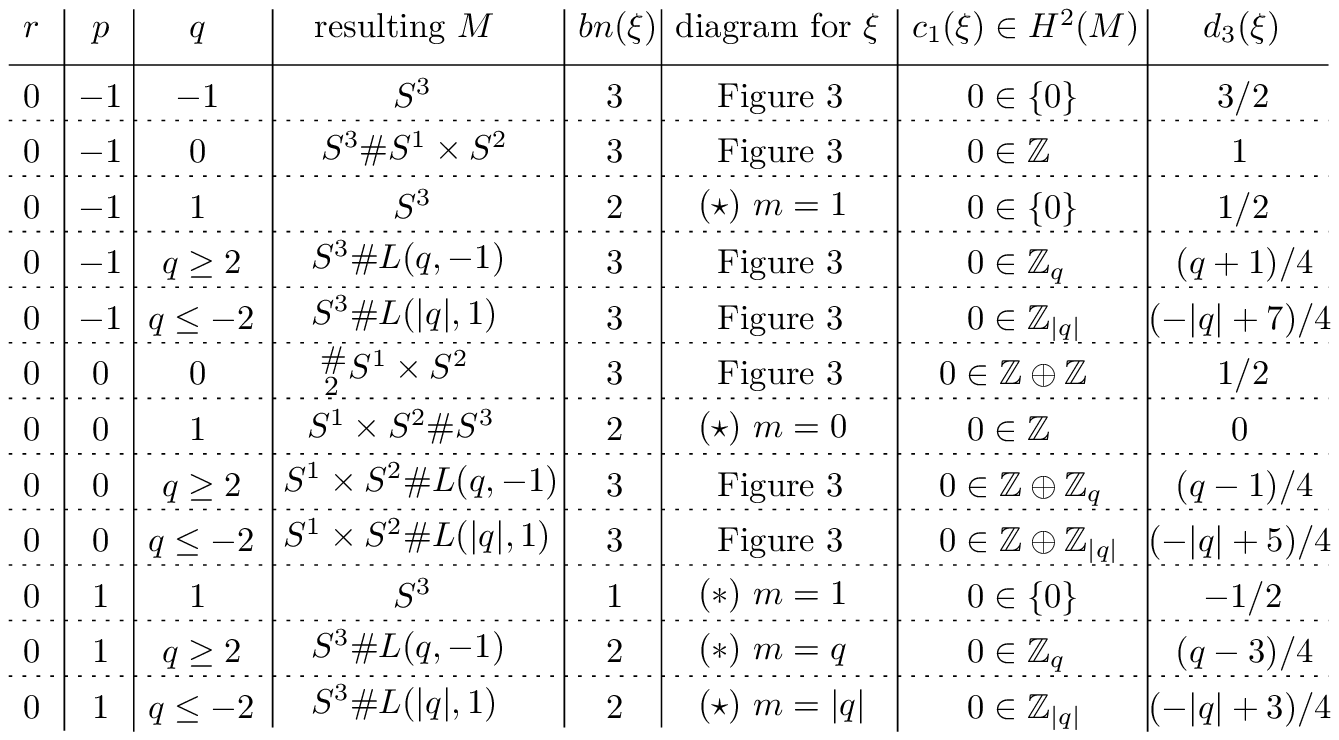}
      \caption{The case $r=0$ ($|p\,|<2$ or $|q|<2$).}
  \label{r_is_0}
\end{table}

\medskip \noindent In \textbf{Table \ref{r_is_0}}, we need to compute
the binding number $bn(\xi)$ for the rows $5$, $12$. For the other rows,
see Remark \ref{Main_Remark}.

\medskip \noindent $\bullet$ If $p=-1, q\leq-2, r=0$, we need to compute
$d_3(\xi_{-1,q,0})$ as $c_1(\xi_{-1,q,0})=0$: We have
$$\mathcal{L}_{-1,q,0}=\left( \begin{array}{cc}
A_1&\mathbf{0}\\\mathbf{0}&A_{|q|+1}\end{array}\right).$$ The
contact structure $\xi_{-1,q,0}$, and $\mathbb{L}_{-1,q,0}$
describing $X_{-1,q,0}$ are given in Figure
\ref{p_-1_q_less_-2_r_0}. We compute that $s=|q|+3$, $c^2=0$,
$\chi(X_{-1,q,0})=|q|+4$, and
$\sigma(X_{-1,q,0})=\sigma(A_1)+\sigma(A_{|q|})=0+|q|-1=|q|-1$, and
so we obtain $d_3(\xi_{-1,q,0})=(-|q|+7)/4$ by Corollary
\ref{d3_computation}. Therefore, $\xi_{-1,q,0}$ is not isotopic to
$\eta_{|q|}$ as $d_3(\eta_{|q|})=(-|q|+3)/4$. Hence,
$bn(\xi_{-1,q,0})=3$ for any $q\leq-2$ by Theorem
\ref{bindingnumber=<2}.

\begin{figure}[ht]
   \includegraphics{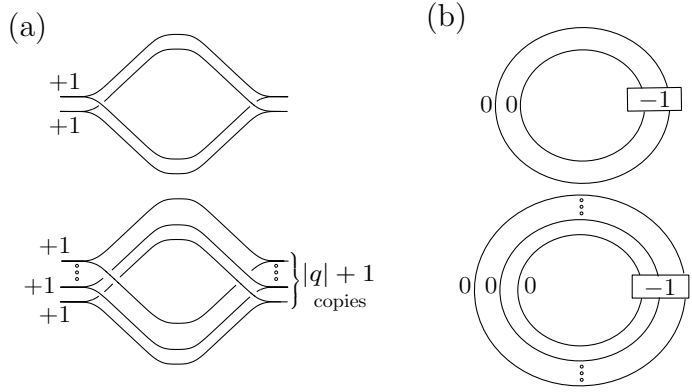}
      \caption{(a) The contact structure $\xi_{-1,q,0}$ on $S^3 \# L(|q|,1) \approx L(|q|,1)$, (b) The corresponding framed link
$\mathbb{L}_{-1,q,0}.$ }
  \label{p_-1_q_less_-2_r_0}
\end{figure}

\medskip \noindent $\bullet$ If $p=1, q\leq-2, r=0$, we have
$(\Sigma,\phi_{1,q,0})=S^{+}_{a}(H^+,D_b^{q})$ (recall the
identification of $\Sigma$ and the curves $a,b,c$ in Figure \ref{our
surface}). Therefore, $\xi_{1,q,0}\cong \eta_{|q|}$ since
$(H^+,D_b^{q})$ supports the overtwisted structure $\eta_{|q|}$ on
$L(|q|,1)$. Hence, $bn(\xi_{1,q,0})=2$ for $q\leq-2$.

\begin{table}[ht]
   \includegraphics{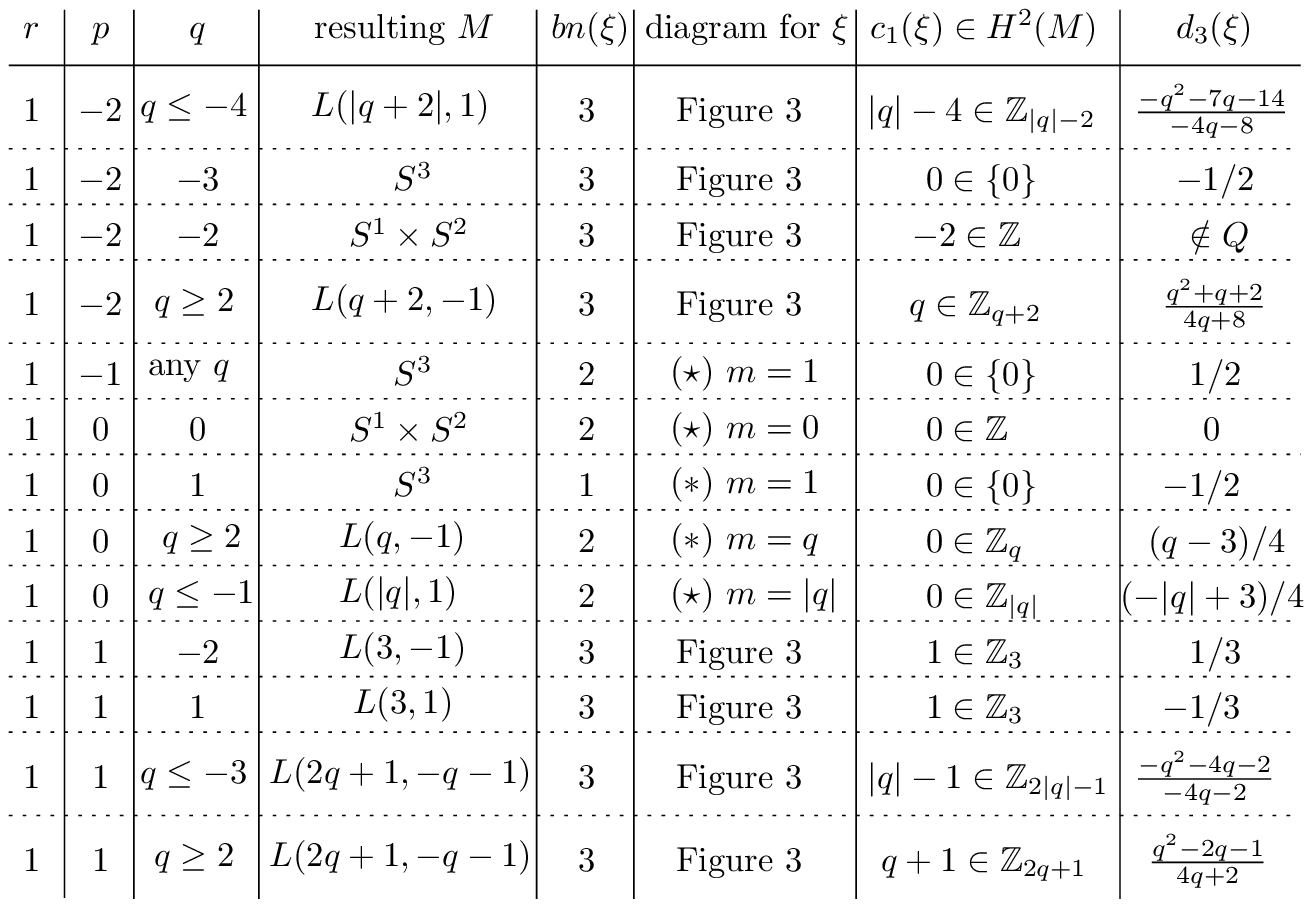}
      \caption{The case $r=1$, $|p\,|<2$ or $|q|<2$ (and the case $(p,q,r)=(-2,q,1)$).}
  \label{r_is_1}
\end{table}

\medskip \noindent In \textbf{Table \ref{r_is_1}}, we need to compute the binding
number $bn(\xi)$ for the rows $1$ and $9$. For the
other rows, see Remark \ref{Main_Remark}.

\medskip \noindent $\bullet$ If $p=-2,q\leq-4, r=1$, let $K_i$'s be the components (with the given orientations) of $\mathbb{L}_{-2,q,1}$ as in Figure \ref{p_-2_q_less_-1_r_1}. Then we obtain the linking matrix
$$\mathcal{L}_{-2,q,1}=\left( \begin{array}{rrrrrrr} -3&-1&-1&-1& -1 &\cdots&-1\\-1& 0 & -1 &-1 & 0&\cdots&0 \\ -1& -1 & 0 &-1 & 0 &\cdots&0 \\-1& -1 & -1 &0 & 0 &\cdots&0 \\-1&0& 0& 0& & & \\ \cdot&\cdot& \cdot &\cdot & & &\\ \cdot&\cdot& \cdot&\cdot & & \Large{A_{|q|}} &  \\ \cdot&\cdot& \cdot& \cdot& & &\\-1&0& 0& 0& & & \\ \end{array}\right)
$$
\begin{figure}[ht]
   \includegraphics{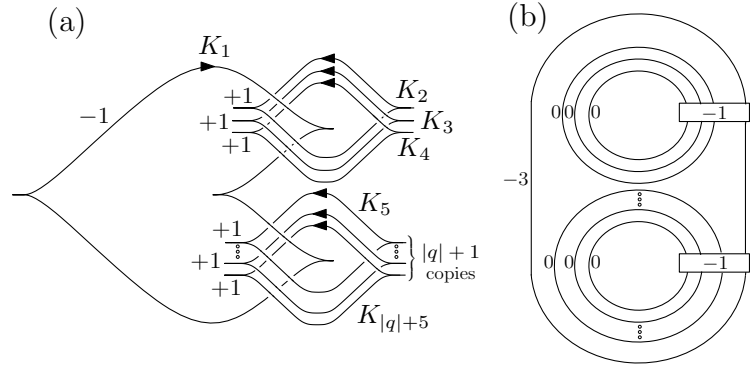}
      \caption{(a) The overtwisted contact structure $\xi_{-2,q,1}$ on $L(|q+2|,1),$ (b) The corresponding framed link $\mathbb{L}_{-2,q,1}.$ }
  \label{p_-2_q_less_-1_r_1}
\end{figure}

\noindent It is not hard to see that $$H_1(M)=\langle \; \mu_1,
\mu_2, \cdots, \mu_{|q|+5} | \;
\mathcal{L}_{-2,q,1}[\mathbf{\mu}]^T_{|q|+5}=[\mathbf{0}]^T_{|q|+5}
\; \rangle=\langle \; \mu_2 | \; (|q|-2)\mu_2=0 \; \rangle \cong
\Z_{|q|-2},$$ and $\mu_1=(|q|-4)\mu_2$. Therefore,
$$c_1(\xi_{-2,q,1})=PD^{-1}(\mu_1)=PD^{-1}(|q|-4)\mu_2=|q|-4 \in
\Z_{|q|-2}.$$ Thus, if $q<-4$, then $\xi_{-2,q,1}$ is not isotopic
to $\eta_{|q+2|}$ as $c_1(\eta_{|q+2|})=0$ implying that
$bn(\xi_{-2,q,1})=3$ by Theorem \ref{bindingnumber=<2}. If $q=-4$,
we compute that $d_3(\xi_{-2,-4,1})=-1/4\neq1/4=d_3(\eta_2)$, so
$bn(\xi_{-2,-4,1})=3$.

\medskip \noindent $\bullet$ If $p=0, q\leq-1, r=1$, we have $(\Sigma,\phi_{0,q,1})=S^{+}_{c}(H^+,D_b^{q})$ (again recall the identification of $\Sigma$ and the curves $a,b,c$ in Figure \ref{our surface}). Therefore, $\xi_{0,q,1}\cong \eta_{|q|}$ since $(H^+,D_b^{q})$ supports the overtwisted structure $\eta_{|q|}$ on $L(|q|,1)$. Hence, $bn(\xi_{0,q,1})=2$ for $q<0$.


\begin{table}[ht]
   \includegraphics{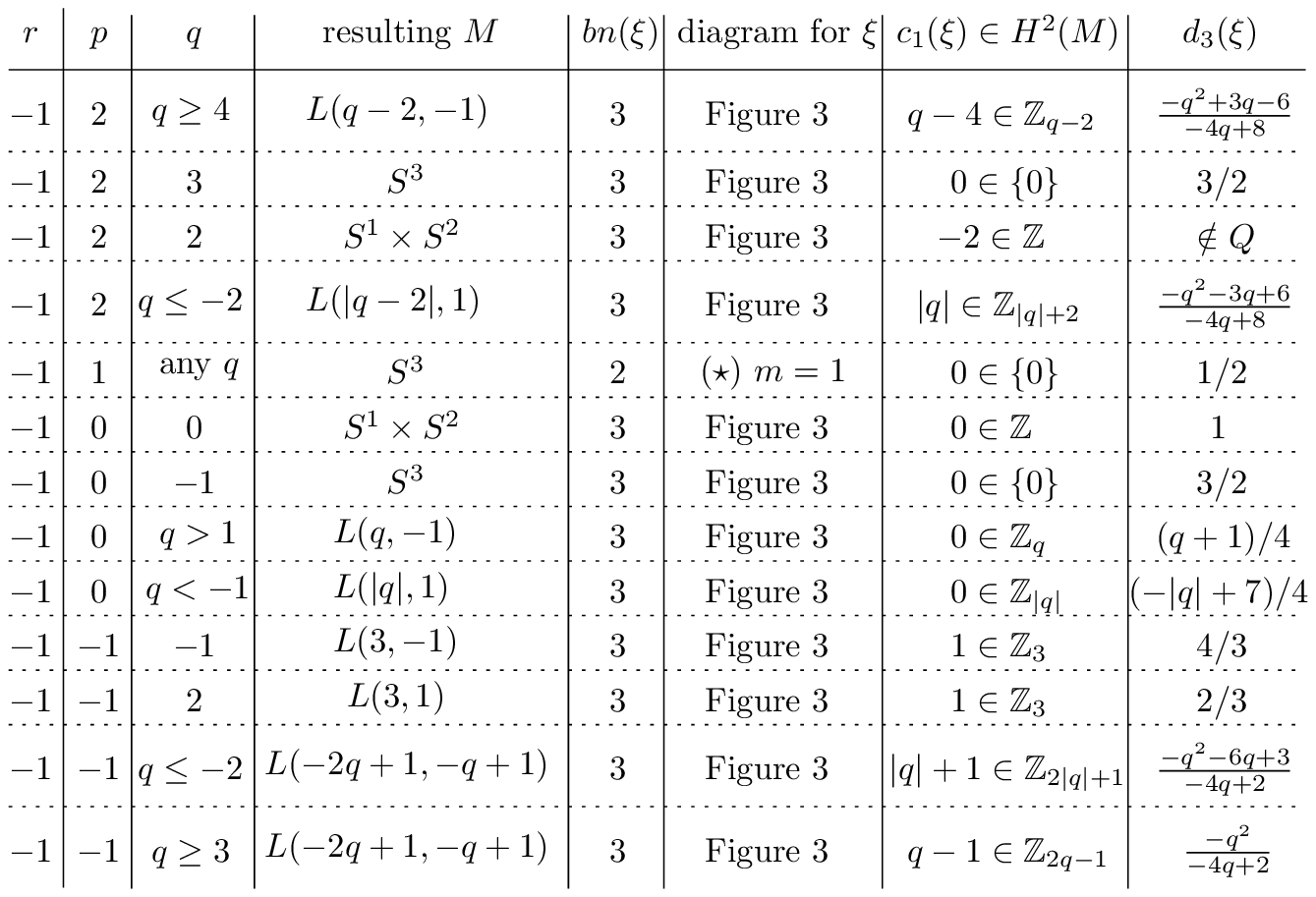}
      \caption{The case $r=-1$, $|p\,|<2$ or $|q|<2$ (and the case $(p,q,r)=(2,q,-1)$)}
  \label{r_is_-1}
\end{table}

\medskip \noindent In \textbf{Table \ref{r_is_-1}}, we need to determine the binding number $bn(\xi)$ for the rows $4$, $7$, $9$, and $11$. For the other rows, see Remark \ref{Main_Remark}.

\medskip \noindent $\bullet$ If $p=2,q\leq-2, r=-1$, then using the corresponding matrix $\mathcal{L}_{2,q,-1}$, we have
$$H_1(M)=\langle \; \mu_1, \mu_2, \cdots, \mu_{|q|+3} | \; \mathcal{L}_{2,q,-1}[\mathbf{\mu}]^T_{|q|+3}=[\mathbf{0}]^T_{|q|+3} \; \rangle=\langle \; \mu_2 | \; (|q|+2)\mu_2=0 \; \rangle \cong \Z_{|q|+2},$$ and $\mu_1=|q|\mu_2$. Therefore, $$c_1(\xi_{2,q,-1})=PD^{-1}(\mu_1)=PD^{-1}(|q|\mu_2)=|q| \in \Z_{|q|+2}.$$ Thus, if $q\leq-2$, then $\xi_{2,q,-1}$ is not isotopic to $\eta_{|q-2|}$ as $c_1(\eta_{|q-2|})=0$ implying that $bn(\xi_{2,q,-1})=3$ by Theorem \ref{bindingnumber=<2}.

\medskip \noindent $\bullet$ If $p=0,q\leq-1, r=-1$ (the rows $7$ or $9$), then $c_1(\xi_{0,q,-1})=0$ and so we need to compute $d_3(\xi_{0,q,-1})$. Let $K_i$'s be the components of $\mathbb{L}_{0,q,-1}$ as in Figure \ref{p_0_q_less_-1_r_-1}. Then
$$\mathcal{L}_{0,q,-1}=\left( \begin{array}{rrrrr} -1&-1&-1&\cdots&-1\\-1&0&0&\cdots&0 \\ -1&0& & & \\ \cdot &\cdot& & & \\ \cdot&\cdot& &\quad \Large{A_{|q|}}&  \\ \cdot&\cdot& & &\\-1&0& & & \\ \end{array}\right)
\longrightarrow \left( \begin{array}{rrrcr}
1&0&0&\cdots&0\\0&-1&0&\cdots&0 \\ 0&0& & & \\ \cdot &\cdot& & & \\
\cdot&\cdot& &\quad A_{|q|}&  \\ \cdot&\cdot& & &\\0&0& & & \\
\end{array}\right)
$$
\noindent By diagonalizing the first two rows of
$\mathcal{L}_{0,q,-1}$, we obtain the matrix on the right. So
$\sigma(\mathcal{L}_{0,q,-1})=\sigma(A_{|q|})=|q|-1$. The contact
surgery diagram for $\xi_{0,q,-1}$ and the corresponding
$4$-manifold $X_{0,q,-1}$ (with $\partial X_{0,q,-1}=M$) are given
in Figure \ref{p_0_q_less_-1_r_-1}.

\begin{figure}[ht]
  \begin{center}
   \includegraphics{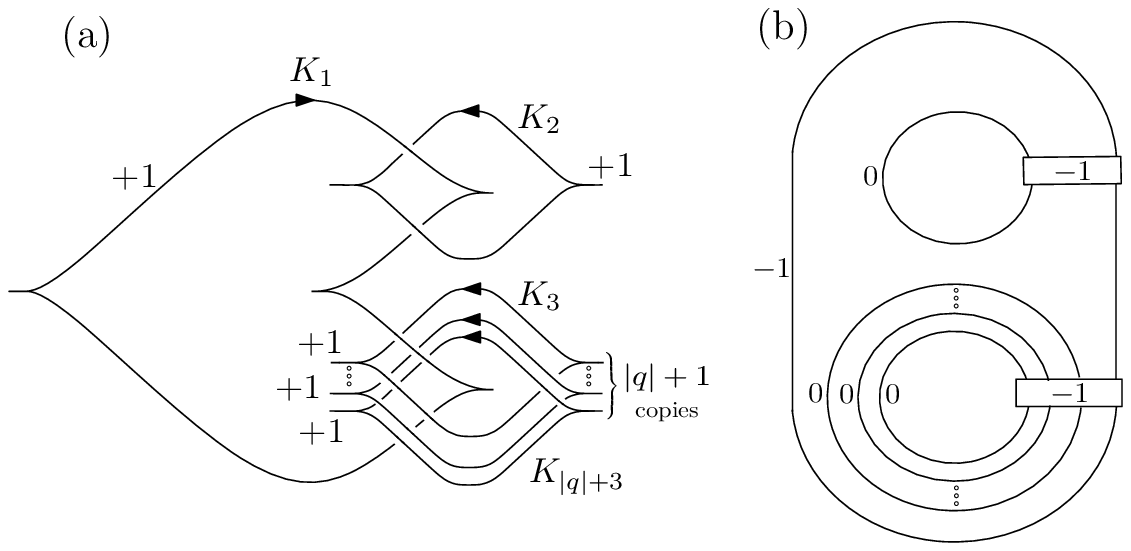}
      \caption{(a) The overtwisted contact structure $\xi_{0,q,-1}$ on $L(|q|,1),$ (b) The corresponding framed link $\mathbb{L}_{0,q,-1}.$ }
  \label{p_0_q_less_-1_r_-1}
    \end{center}
\end{figure}

\noindent Then the system
$\mathcal{L}_{0,q,-1}[\mathbf{b}]^T=[rot(K_1) \;rot(K_2) \cdots
rot(K_{|q|+3})]^T=[1 \; 0 \; 0 \cdots 0]^T$ has the solution
$[\mathbf{b}]=[0 \; -1 \;  0 \cdots 0]$, and so $c^2=0$. Moreover,
$\chi(X_{0,q,-1})=|q|+4$ and $s=|q|+3$. Therefore, we obtain
$d_3(\xi_{0,q,-1})=(-|q|+7)/4$ implying that $\xi_{0,q,-1}$ is not
isotopic to $\eta_{|q|}$ as $d_3(\eta_{|q|})=(-|q|+3)/4$. Hence,
$bn(\xi_{0,q,-1})=3$ by Theorem \ref{bindingnumber=<2}.

\medskip \noindent $\bullet$ If $p=-1,q=2, r=-1$, we have $c_1(\xi_{-1,2,-1})=1$ implying that $bn(\xi_{-1,2,-1})=3$. To see this, note that $c_1(\xi{-1,2,-1})=PD^{-1}(\mu_1)$ where $\mu_1$ is the meridian of the surgery curve corresponding $K_1$. Then using
$$\mathcal{L}_{-1,2,-1}=\left( \begin{array}{rrrr} -1&-1&-1&-1\\-1&-2&0&0 \\ -1&0&0&-1 \\ -1&0&-1&0\\ \end{array}\right)$$
we get $H_1(M)=\langle \; \mu_1, \mu_2, \mu_3, \mu_4 | \;
\mathcal{L}_{-1,2,-1}[\mathbf{\mu}]^T_{4}=[\mathbf{0}]^T_{4} \;
\rangle=\langle \; \mu_2 | \; 3\mu_2=0 \; \rangle \cong \Z_{3},$ and
$\mu_1=-2\mu_2$. Therefore, we compute
$$c_1(\xi_{-1,2,-1})=PD^{-1}(\mu_1)=PD^{-1}(-2\mu_2)=-2 \in
\Z_{3}\equiv 1 \in \Z_{3}.$$


\begin{table}[ht]
   \includegraphics{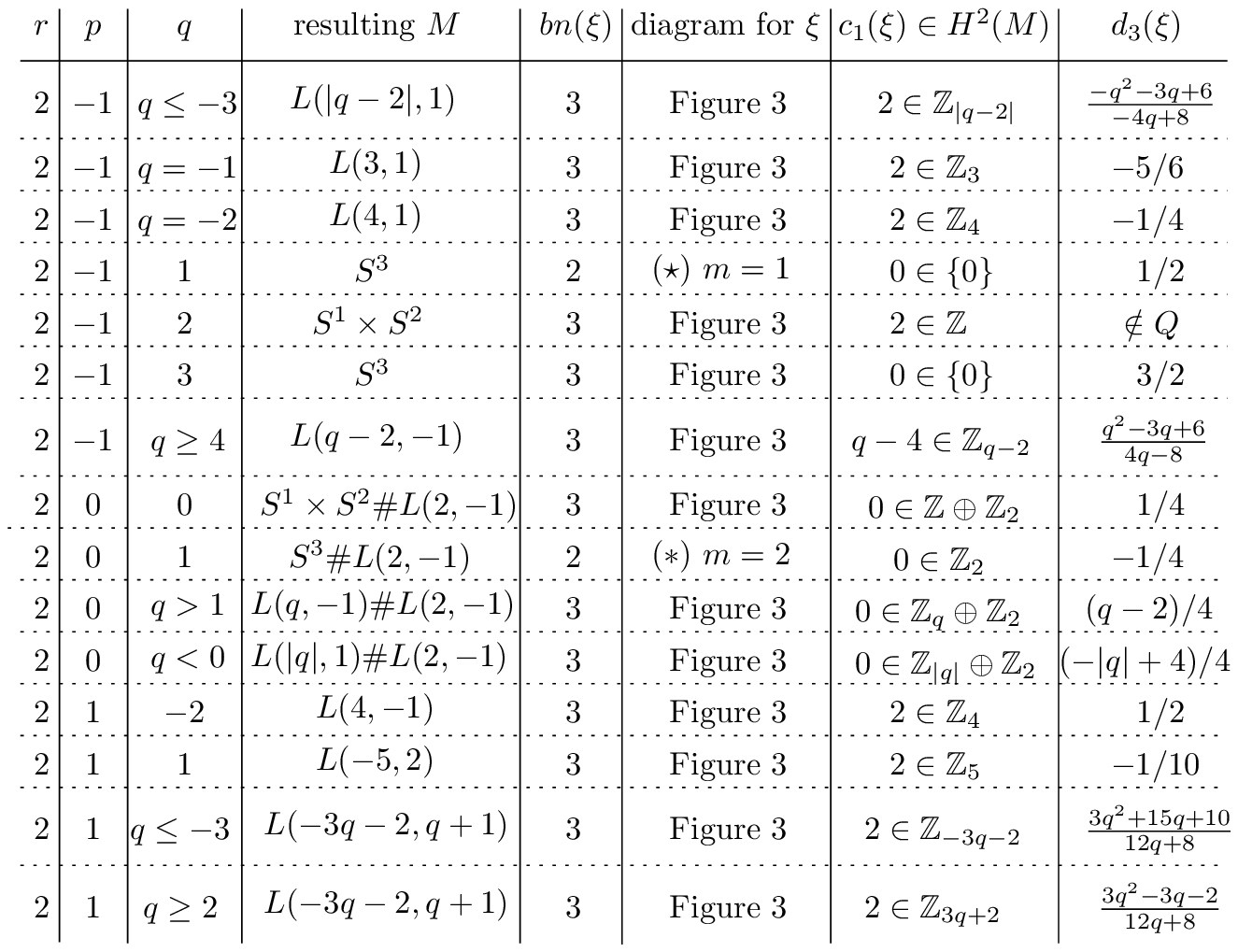}
      \caption{The case $r=2$ ($|p\,|<2$ or $|q|<2$).}
  \label{r_is_2}
\end{table}

\medskip \noindent In \textbf{Table \ref{r_is_2}}, we need to compute the binding number
$bn(\xi)$ for the rows $1$, $2$, and $3$. For the other rows, again
see Remark \ref{Main_Remark}.

\medskip \noindent For the first three rows in Table \ref{r_is_2}, the contact structure $\xi_{-1,q,2}$ on $L(|q-2|,1)$ and the link $\mathbb{L}_{-1,q,2}$ ($q\leq -1$) are given in Figure \ref{p_-1_q_less_-1_r_2}. We write the linking matrix $\mathcal{L}_{-1,q,2}$ as the matrix on the left below. It is not hard to see that $c_1(\xi_{-1,q,2})=2 \in \Z_{|q-2|}$, and so $bn(\xi_{-1,q,2})=3$. As an illustration we will compute $d_3(\xi_{-1,q,2})$ (even though it is not necessary for the proof). The matrix on the right below is obtained by diagonalizing the first two rows of $\mathcal{L}_{-1,q,2}$. So we compute $\sigma(\mathcal{L}_{-1,q,2})=2+\sigma(A_1)+\sigma(B_{|q|})$ which is $|q|-1$ if $q\leq-3$, and is equal to 2 if $q=-1,-2$ (recall $\sigma(B_n)$ is $n-3$ if $n\geq3$, and 0 if $n=1,2$).

\begin{figure}[ht]
   \includegraphics{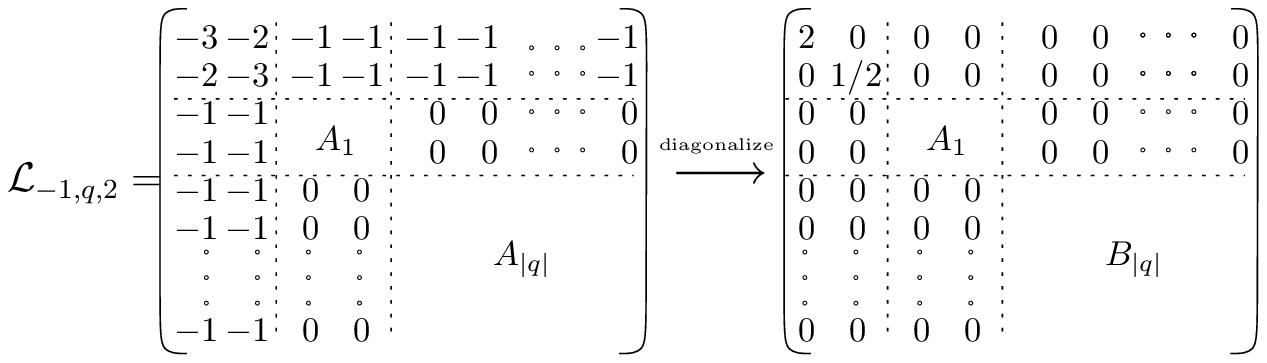}
  \label{matrices_2}
\end{figure}

\begin{figure}[ht]
   \includegraphics{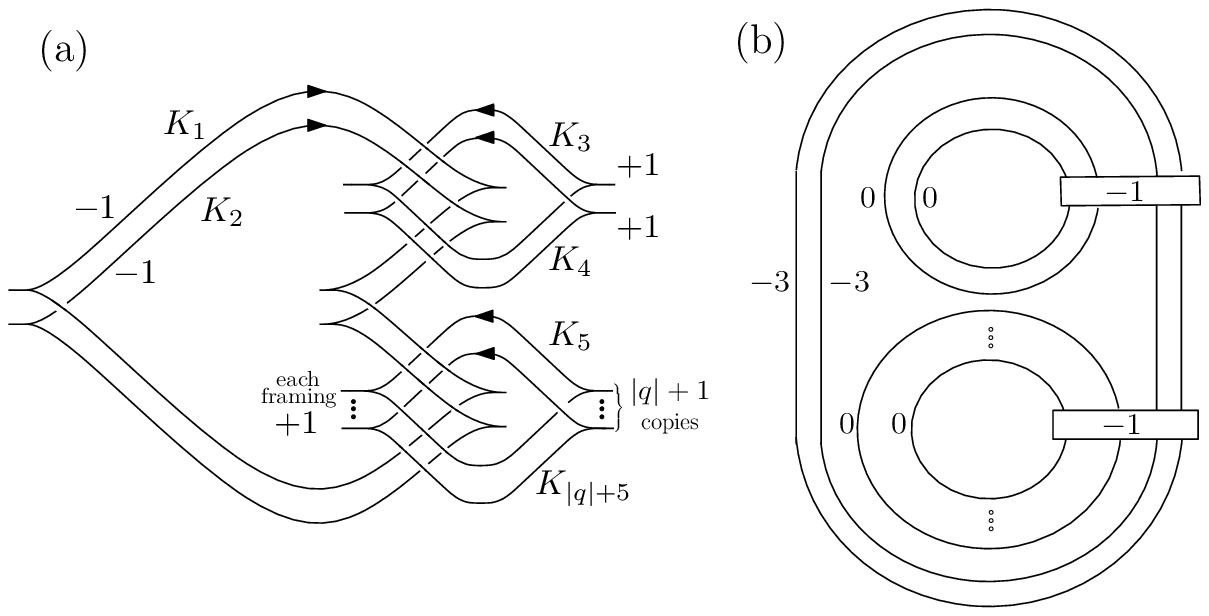}
      \caption{(a) The contact structure $\xi_{-1,q,2}$ on $L(|q-2|,1)$ for $q\leq-1$, (b) The corresponding framed link $\mathbb{L}_{-1,q,2}.$ }
  \label{p_-1_q_less_-1_r_2}
\end{figure}

\noindent By a standard calculation, the system
$$\mathcal{L}_{-1,q,2}[\mathbf{b}]^T=[rot(K_1) \;rot(K_2) \cdots
rot(K_{|q|+5})]^T=[1 \; 1 \; 0 \cdots 0]^T$$ has the solution
$[\mathbf{b}]=[\frac{|q|}{|q|+2} \; \frac{|q|}{|q|+2} \;
\frac{-2|q|}{|q|+2} \; \frac{-2|q|}{|q|+2} \; \frac{-2}{|q|+2}
\cdots \frac{-2}{|q|+2}]$ for $q\leq-1$, and so we compute
$$c^2=[\mathbf{b}]\mathcal{L}_{-1,q,2}[\mathbf{b}]^T=2|q|/(|q|+2).$$

\medskip \noindent $\bullet$ If $p=-1,q=-1, r=2$, then $c^2=2/3$, $\sigma(X_{-1,-1,2})=2$,
$\chi(X_{-1,-1,2})=7$, and $s=4$. So we get
$d_3(\xi_{-1,-1,2})=-5/6$.

\medskip \noindent $\bullet$ If $p=-1,q=-2, r=2$, then $c^2=1$, $\sigma(X_{-1,-2,2})=2$,
$\chi(X_{-1,-2,2})=8$, and $s=5$. Therefore, we get
$d_3(\xi_{-1,-2,2})=-1/4$.

\medskip \noindent $\bullet$ If $p=-1,q\leq-3, r=2$, then $c^2=2|q|/(|q|+2)$, $\sigma(X_{-1,q,2})=|q|-1$,
$\chi(X_{-1,q,2})=|q|+6$, and $s=|q|+3$. So we obtain
$$d_3(\xi_{-1,q,2})=\displaystyle{\frac{-q^2-3q+6}{-4q+8}}.$$

\begin{table}[ht]
  \begin{center}
   \includegraphics{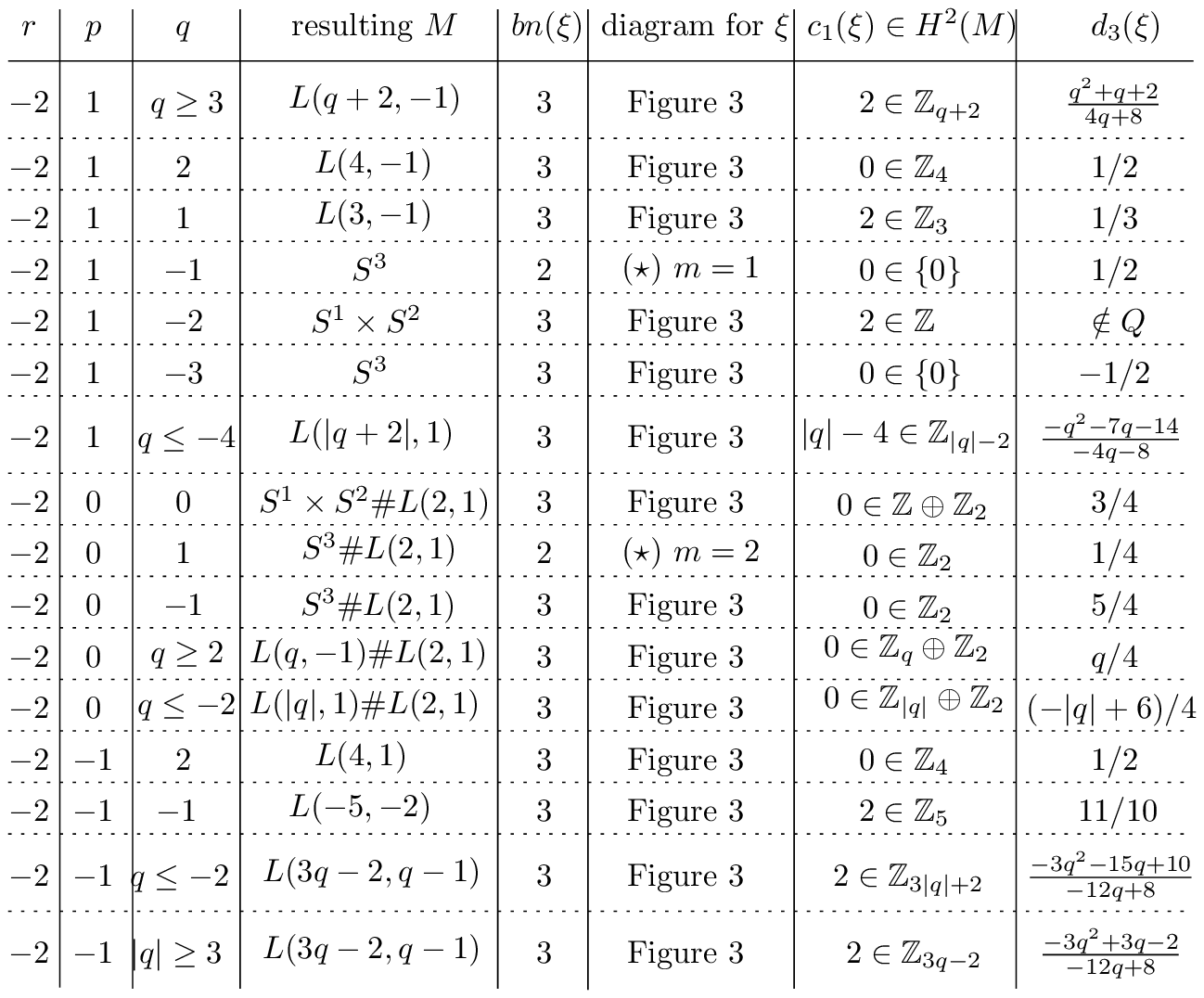}
      \caption{The case $r=-2$ ($|p\,|<2$ or $|q|<2$).}
  \label{r_is_-2}
    \end{center}
\end{table}

\medskip \noindent In \textbf{Table \ref{r_is_-2}}, we need to compute the binding number $bn(\xi)$ for the rows $7$, $9$, $10$, and $13$. For the other rows, see Remark \ref{Main_Remark}.

\medskip \noindent $\bullet$ If $p=1,q\leq-4, r=-2$, the contact structure $\xi_{1,q,-2}$ on $L(|q+2|,1)$ and the link $\mathbb{L}_{1,q,-2}$ are given in Figure \ref{p_1_q_less_-3_r_-2}.

\begin{figure}[ht]
  \begin{center}
   \includegraphics{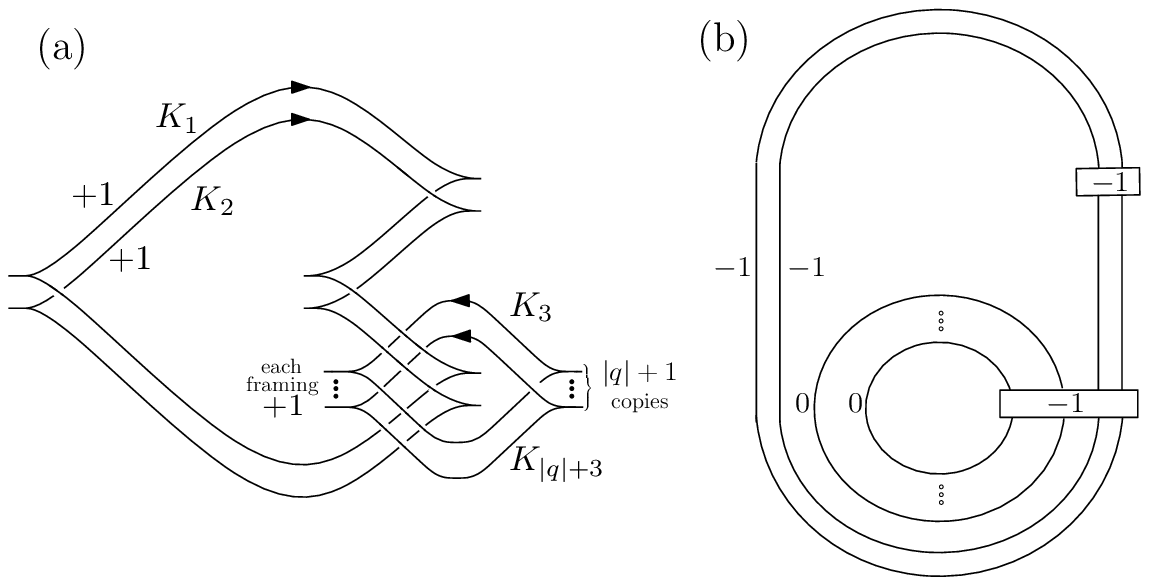}
      \caption{(a) The contact structure $\xi_{1,q,-2}$ on $L(|q+2|,1)$ for $q<-3$,  (b) The corresponding framed link $\mathbb{L}_{1,q,-2}.$ }
  \label{p_1_q_less_-3_r_-2}
    \end{center}
\end{figure}

\noindent We will first compute that $c_1(\xi_{1,q,-2})=|q|-4 \in \Z_{|q|-2}$ (so $bn(\xi_{1,q,-2})=3$), and then (even though it is not necessary for the proof) we will evaluate $d_3(\xi_{1,q,-2})$ as an another sample computation. Using $\mathcal{L}_{1,q,-2}$ (on the left below), we have
$$\begin{array}{lll}
H_1(M)&=&\langle \; \mu_1, \mu_2, \cdots, \mu_{|q|+3} | \; \mathcal{L}_{1,q,-2}[\mathbf{\mu}]^T_{|q|+3}=[\mathbf{0}]^T_{|q|+3} \; \rangle \\
      &=&\langle \; \mu_1, \mu_3 | \; -3\mu_1-(|q|+1)\mu_3=0, -2\mu_1-|q|\mu_3=0 \; \rangle \\
      &=&\langle \; \mu_3 | \; (|q|-2)\mu_3=0 \; \rangle \cong \Z_{|q|-2},
\end{array}$$
and also we have $\mu_1=\mu_2=-\mu_3$. Therefore, we obtain
$$c_1(\xi_{2,q,-1})=PD^{-1}(\mu_1+\mu_2)=PD^{-1}(-2\mu_3)=-2\equiv
|q|-4 \in \Z_{|q|-2}.$$

\noindent The matrix on the right below is obtained by diagonalizing
the first two rows of $\mathcal{L}_{1,q,-2}$. So we compute
$\sigma(\mathcal{L}_{1,q,-2})=0+ \sigma(C_{|q|})=|q|-1$ (recall
$\sigma(C_n)=n-1$ if $n\geq2$).

\begin{figure}[ht]
  \begin{center}
   \includegraphics{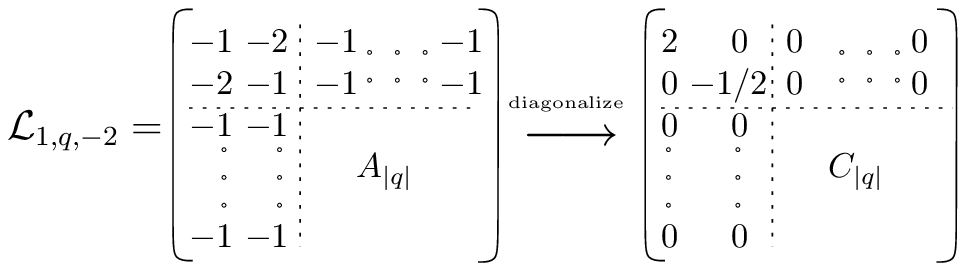}
  \label{matrices_3}
    \end{center}
\end{figure}

\noindent By a standard calculation, the system
$$\mathcal{L}_{1,q,-2}[\mathbf{b}]^T=[rot(K_1) \;rot(K_2) \cdots
rot(K_{|q|+3})]^T=[1 \; 1 \; 0 \cdots 0]^T$$ has the solution
$[\mathbf{b}]=[\frac{-|q|}{|q|-2} \; \frac{-|q|}{|q|-2} \;
\frac{2}{|q|-2} \; \cdots \frac{2}{|q|-2}]$, and so we obtain
$$c^2=[\mathbf{b}]\mathcal{L}_{1,q,-2}[\mathbf{b}]^T=-2|q|/(|q|-2).$$

\medskip \noindent Moreover, $\chi(X_{1,q,-2})=|q|+4$, and $s=|q|+3$. So we compute
$$d_3(\xi_{1,q,-2})=\displaystyle{\frac{-q^2-7q-14}{-4q-8}}.$$

\medskip \noindent $\bullet$ If $p=0,q=1, r=-2$, then $\xi_{0,1,-2}$ and $\mathbb{L}_{0,1,-2}$ are given in Figure \ref{p_0_q_1_r_-2}.

\begin{figure}[ht]
  \begin{center}
   \includegraphics{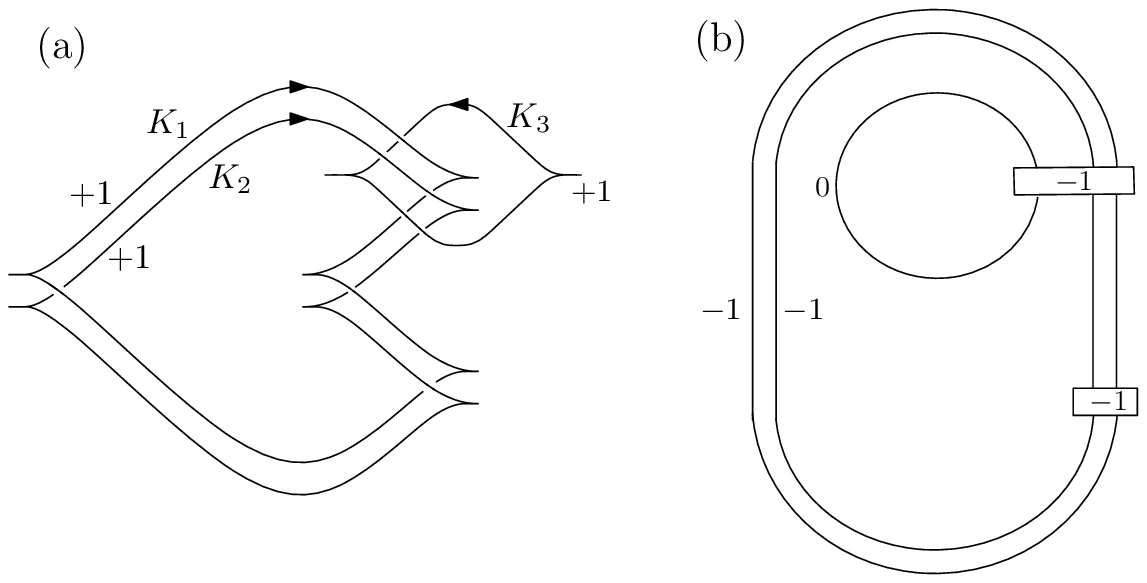}
      \caption{(a) The contact structure $\xi_{0,1,-2}$ on $S^3 \# L(2,1)\approx L(2,1),$ (b) The corresponding framed link $\mathbb{L}_{0,1,-2}.$ }
  \label{p_0_q_1_r_-2}
    \end{center}
\end{figure}

\noindent One can get $c_1(\xi_{0,1,-2})=0$, so we need
$d_3(\xi_{0,1,-2})$. The corresponding linking matrix is
$$\mathcal{L}_{0,1,-2}=\left( \begin{array}{rrr} -1&-2&-1\\-2&-1&-1\\-1&-1&0\\ \end{array}\right)
\longrightarrow \left( \begin{array}{rrr} -1&0&0\\0&2&0\\0&0&1\\
\end{array}\right).
$$
\noindent We diagonalize $\mathcal{L}_{0,1,-2}$, and obtain the
matrix on the right. So $\sigma(\mathcal{L}_{0,1,-2})=1$. We find
that the system $\mathcal{L}_{0,1,-2}[\mathbf{b}]^T=[rot(K_1)
\;rot(K_2) \; rot(K_3)]^T=[1 \; 1 \; 0]^T$ has the solution
$[\mathbf{b}]=[0 \; 0 \; -1]$, and so $c^2=0$. Also we have
$\chi(X_{0,1,-2})=4$ and $s=3$. So we get
$d_3(\xi_{0,1,-2})=1/4=d_3(\eta_{2})$ which implies that
$\xi_{0,1,-2}$ is isotopic to $\eta_{2}$. Thus, $bn(\xi_{0,1,-2})=2$
by Theorem \ref{bindingnumber=<2}.

\medskip \noindent $\bullet$ If $p=0,q=-1, r=-2$, then the contact structure $\xi_{0,1,-2}$ on $L(2,1)$ and the link $\mathbb{L}_{0,-1,-2}$ describing $X_{0,-1,-2}$ are given in Figure \ref{p_0_q_-1_r_-2}. It is easy to check  $c_1(\xi_{0,-1,-2})=0$, so we compute $d_3(\xi_{0,-1,-2})$:

\begin{figure}[ht]
  \begin{center}
   \includegraphics{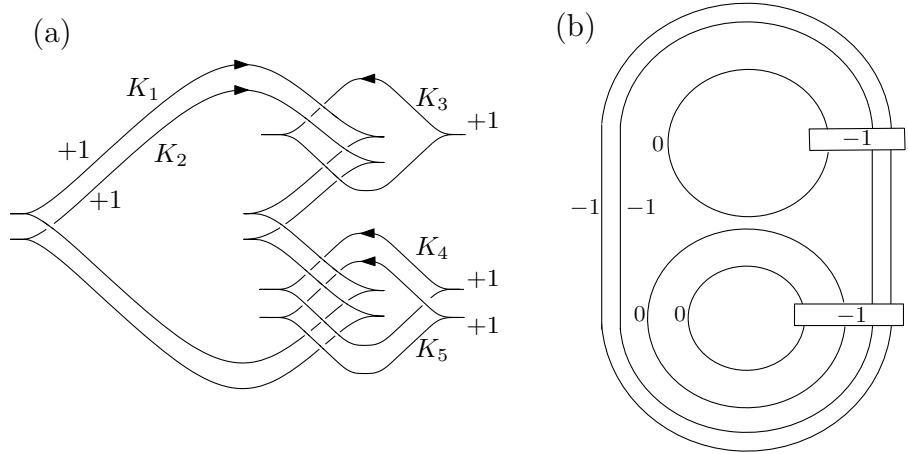}
      \caption{(a) The contact structure $\xi_{0,-1,-2}$ on $S^3 \# L(2,1)\approx L(2,1),$ (b) The corresponding framed link $\mathbb{L}_{0,-1,-2}.$ }
  \label{p_0_q_-1_r_-2}
    \end{center}
\end{figure}

\noindent The corresponding linking matrix is
$$\mathcal{L}_{0,-1,-2}=\left( \begin{array}{rrrrr} -1&-2&-1&-1&-1\\-2&-1&-1&-1&-1\\-1&-1&0&0&0\\-1&-1&0&0&-1\\-1&-1&0&-1&0\\ \end{array}\right)
\longrightarrow \left( \begin{array}{crrcc}
2&0&0&0&0\\0&-1&0&0&0\\0&0& \;\;1&0&0\\0&0&0&-1/2&0\\0&0&0&0&2\\
\end{array}\right).
$$
\noindent We diagonalize $\mathcal{L}_{0,-1,-2}$, and obtain the
matrix on the right. So $\sigma(\mathcal{L}_{0,-1,-2})=1$. The
system $$\mathcal{L}_{0,-1,-2}[\mathbf{b}]^T=[rot(K_1) \;rot(K_2) \;
rot(K_3)\;rot(K_4)\;rot(K_5)]^T=[1 \; \; 1 \; \; 0 \; \; 0 \;
\;0]^T$$ has the solution $[\mathbf{b}]=[0 \; \; 0 \; -1 \; \; 0 \;
\; 0]$ which yields $c^2=0$. Also we have $\chi(X_{0,-1,-2})=4$ and
$s=3$. So we get $d_3(\xi_{0,-1,-2})=5/4\neq1/4=d_3(\eta_{2})$.
Therefore, $\xi_{0,-1,-2}$ is not isotopic to $\eta_{2}$, and so
$bn(\xi_{0,-1,-2})=3$ by Theorem \ref{bindingnumber=<2}.

\medskip \noindent $\bullet$ If $p=-1,q=2, r=-2$, then the contact structure $\xi_{-1,2,-2}$ on $L(4,1)$ and the link $\mathbb{L}_{-1,2,-2}$ describing $X_{-1,2,-2}$ are given in Figure \ref{p_-1_q_2_r_-2}. We compute that $c_1(\xi_{-1,2,-2})=0$, so we need to find $d_3(\xi_{-1,2,-2})$.

\begin{figure}[ht]
  \begin{center}
   \includegraphics{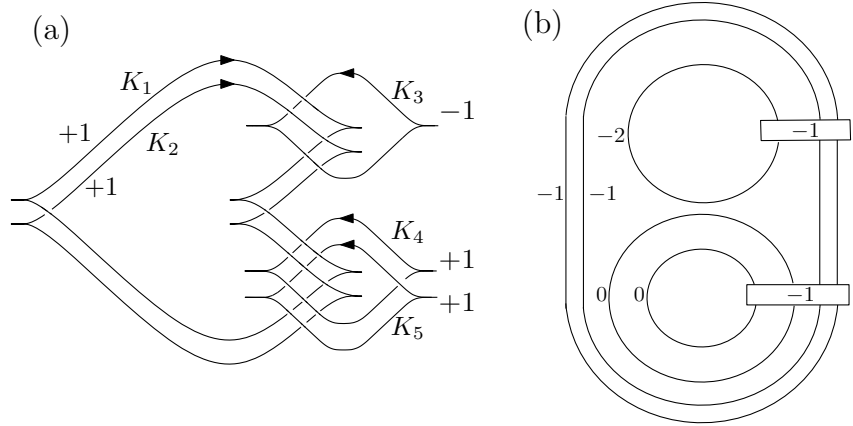}
      \caption{(a) The contact structure $\xi_{-1,2,-2}$ on $L(4,1), \quad \quad \quad \quad \quad \quad$
      (b) The corresponding framed link $\mathbb{L}_{-1,2,-2}.$ }
  \label{p_-1_q_2_r_-2}
    \end{center}
\end{figure}

\noindent The corresponding linking matrix is
$$\mathcal{L}_{-1,2,-2}=\left( \begin{array}{rrrrr} -1&-2&-1&-1&-1\\-2&-1&-1&-1&-1\\-1&-1&-2&0&0\\-1&-1&0&0&-1\\-1&-1&0&-1&0\\ \end{array}\right)
\longrightarrow \left( \begin{array}{crrcc}
2&0&0&0&0\\0&-1&0&0&0\\0&0& -1&0&0\\0&0&0&1&0\\0&0&0&0&2\\
\end{array}\right).
$$
\noindent We diagonalize $\mathcal{L}_{-1,2,-2}$, and obtain the
matrix on the right. So $\sigma(\mathcal{L}_{-1,2,-2})=1$. The
system $$\mathcal{L}_{-1,2,-2}[\mathbf{b}]^T=[rot(K_1) \;rot(K_2) \;
rot(K_3)\;rot(K_4)\;rot(K_5)]^T=[1 \; \; 1 \; \; 0 \; \;0 \; \;0]^T$$
has the solution $[\mathbf{b}]=[1/2 \; \; 1/2 \; \; -1/2 \; \; -1 \;
\; -1]$, so we compute $c^2=1$. Moreover, $\chi(X_{-1,2,-2})=6$ and
$s=4$. Then we get $d_3(\xi_{-1,2,-2})=1/2\neq-1/4=d_3(\eta_{4})$.
Therefore, $\xi_{-1,2,-2}$ is not isotopic to $\eta_{4}$, and so
$bn(\xi_{-1,2,-2})=3$ by Theorem \ref{bindingnumber=<2}.

\begin{table}[ht]
  \begin{center}
   \includegraphics{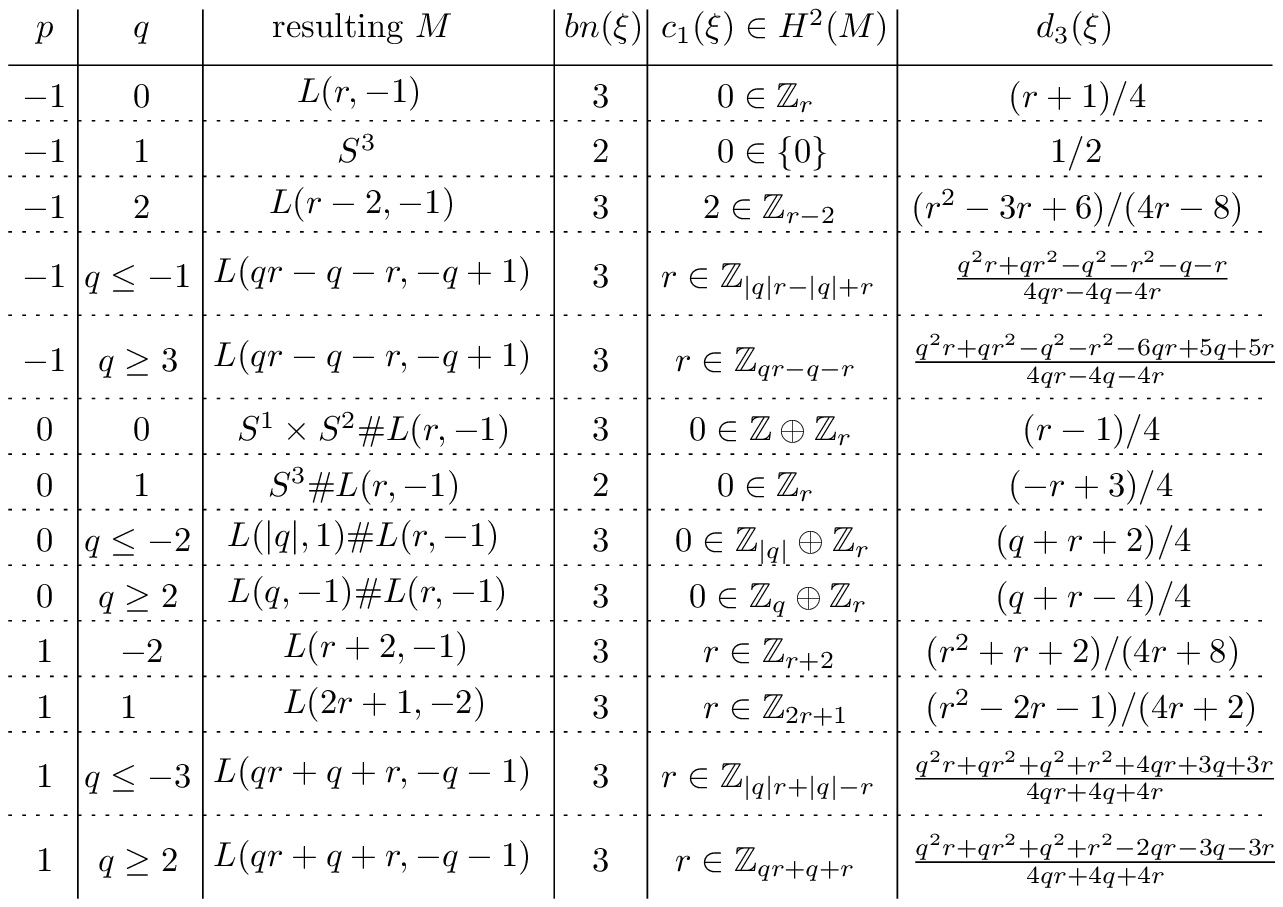}
      \caption{The case $r>2$ ($|p\,|<2$ or $|q|<2$).}
  \label{r_is_greater_2}
    \end{center}
\end{table}

\medskip \noindent In \textbf{Table \ref{r_is_greater_2}}, we do not need any computation to find $bn(\xi)$: For any row, we can use Remark \ref{Main_Remark}. For example, in the $1^{st}$ row, we have an overtwisted contact structure on the lens space $L(m,-1)$ for some $m\geq1$. Therefore, the resulting contact manifold is not listed in Theorem \ref{bindingnumber=<2}, and hence we must have $bn(\xi)=3$.

\begin{table}[ht]
   \includegraphics{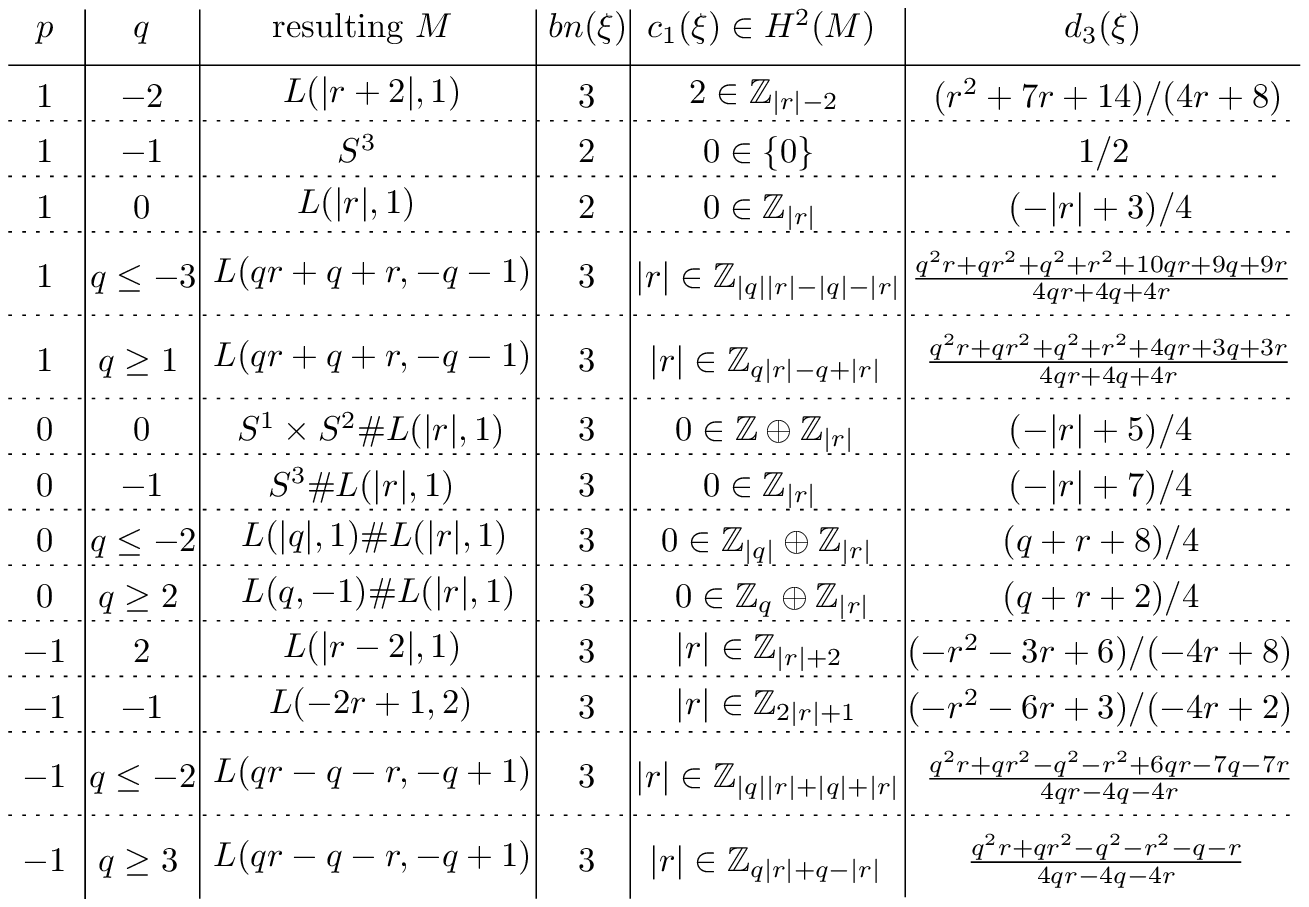}
      \caption{The case $r<-2$ ($|p\,|<2$ or $|q|<2$).}
  \label{r_is_less_-2}
\end{table}

\medskip \noindent In \textbf{Table \ref{r_is_less_-2}}, we need to compute the binding number $bn(\xi)$ for the rows $1$, $3$, $7$, and $10$. For the other rows, see Remark \ref{Main_Remark}.

\medskip \noindent $\bullet$ If $p=1,q=-2, r<-2$, $\xi_{1,-2,r}$ is an overtwisted contact structure  on $L(|r+2|,1)$. It is not hard to see that
$c_1(\xi_{1,-2,r})=2 \in \Z_{|r|-2}$. Therefore, we immediately get $bn(\xi_{1,q,-2})=3$ because $c_1(\eta_{|r+2|})=0$.

\medskip \noindent $\bullet$ If $p=1,q=0, r<-2$, the contact structure $\xi_{1,0,r}$
on $L(|r|,1)$ and the link $\mathbb{L}_{1,0,r}$ are given in Figure
\ref{p_1_q_0_r_less_-2}. It is easy to see that $c_1(\xi_{1,0,r})=0
\in \Z_{|r|}$, so we need $d_3(\xi_{1,0,r})$: The corresponding linking matrix is on the left below.
Diagonalize $\mathcal{L}_{1,0,r}$ to get the matrix on the right.
Therefore, $\sigma(\mathcal{L}_{1,0,r})=|r|-1$.

\begin{figure}[ht]
   \includegraphics{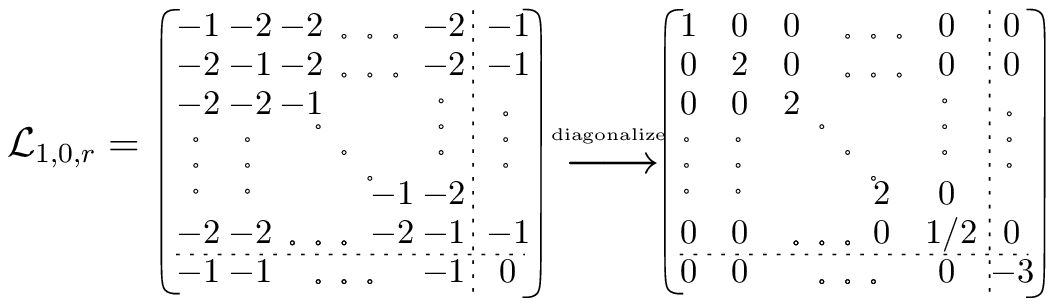}
  \label{matrices_5}
\end{figure}

\begin{figure}[ht]
   \includegraphics{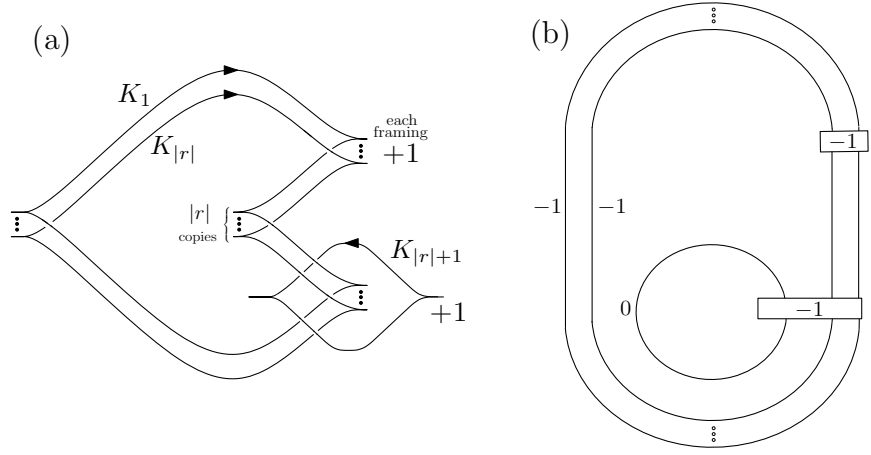}
      \caption{(a) The contact structure $\xi_{1,0,r}$ on $L(|r|,1)$ for $r<-2, \quad$ (b) The corresponding framed link $\mathbb{L}_{1,0,r}.$ }
  \label{p_1_q_0_r_less_-2}
\end{figure}

\noindent The system $$\mathcal{L}_{1,0,r}[\mathbf{b}]^T=[rot(K_1)
\;rot(K_2) \cdots rot(K_{|q|+1})]^T=[1 \cdots 1 \; 0]^T$$ has the
solution $[\mathbf{b}]=[0 \cdots 0 \; -1]$, and so we obtain
$c^2=0$. Moreover, $\chi(X_{1,0,r})=|r|+2$, and $s=|r|+1$. So we
compute $d_3(\xi_{1,0,r})=(-|r|+3)/4=d_3(\eta_{|r|})$ which implies
that $\xi_{1,0,r}$ is isotopic to $\eta_{|r|}$ on $L(|r|,1)$. Thus,
$bn(\xi_{1,0,r})=2$ by Theorem \ref{bindingnumber=<2}.

\medskip \noindent $\bullet$ If $p=0,q=-1, r<-2$, the contact structure $\xi_{0,-1,r}$ on $L(|r|,1)$ and the link $\mathbb{L}_{0,-1,r}$ are given in Figure \ref{p_0_q_-1_r_less_-2}. Again we have $c_1(\xi_{0,-1,r})=0 \in \Z_{|r|}$, so we need to find $d_3(\xi_{0,-1,r})$: We diagonalize $\mathcal{L}_{0,-1,r}$ and get the matrix on the right below. So, we conclude that $\sigma(\mathcal{L}_{0,-1,r})=|r|-1$.

\begin{figure}[ht]
   \includegraphics{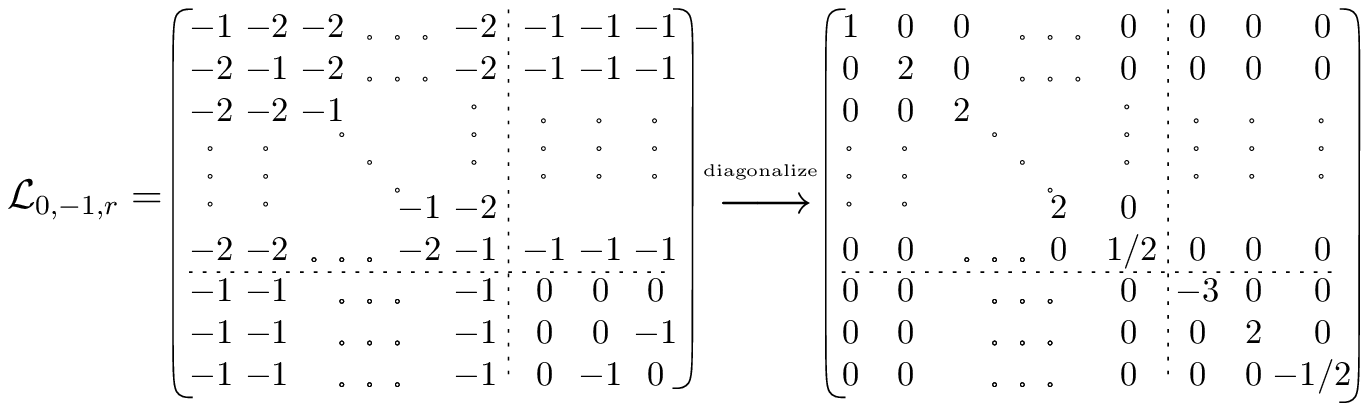}
  \label{matrices_6}
\end{figure}

\begin{figure}[ht]
  \begin{center}
   \includegraphics{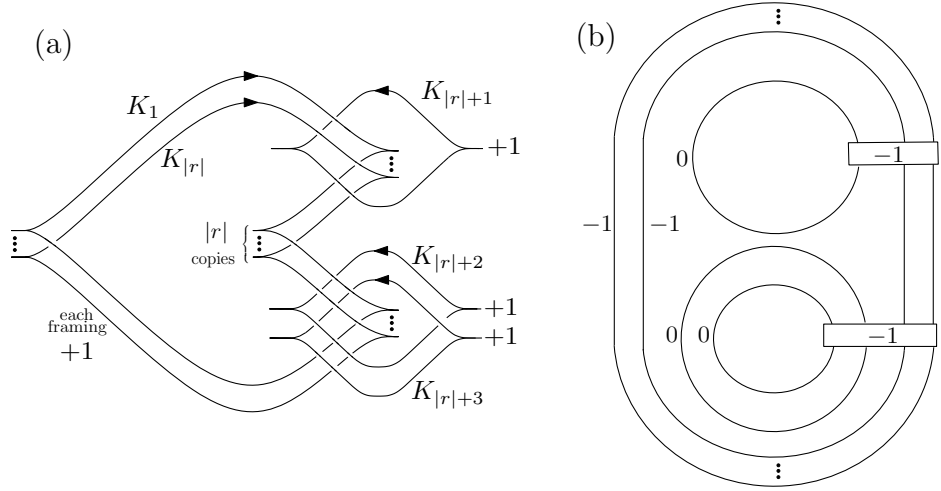}
      \caption{(a) The contact structure $\xi_{0,-1,r}$ on $L(|r|,1)$ for $r<-2,\quad \quad$  (b) The corresponding framed link $\mathbb{L}_{0,-1,r}.$ }
  \label{p_0_q_-1_r_less_-2}
    \end{center}
\end{figure}

\noindent The system $$\mathcal{L}_{0,-1,r}[\mathbf{b}]^T=[rot(K_1)
\;rot(K_2) \cdots rot(K_{|q|+3})]^T=[1 \cdots 1 \; 0 \; 0 \; 0]^T$$
has the solution $[\mathbf{b}]=[0 \cdots \;0 \; 1 \; 0 \;0]$, so we
get $c^2=0$. Also $\chi(X_{1,q,-2})=|r|+4$, and $s=|r|+3$. So we
compute $d_3(\xi_{0,-1,r})=(-|r|+7)/4 \neq
(-|r|+3)/4=d_3(\eta_{|r|})$ implying that $\xi_{0,-1,r} \ncong \eta_{|r|}$
on $L(|r|,1)$. Hence, $bn(\xi_{0,-1,r})=3$ by Theorem \ref{bindingnumber=<2}.

\medskip \noindent $\bullet$ If $p=-1,q=2, r<-2$, we have $bn(\xi_{-1,2,r})=3$
because $c_1(\xi_{-1,2,r})=|r| \in \Z_{|r|+2}$. We compute
$c_1(\xi_{-1,2,r})$ as follows: We use the linking matrix
$\mathcal{L}_{-1,2,r}$ to get the representation
$$\begin{array}{lll}
H_1(M)&=&\langle \; \mu_1, \mu_2, \cdots, \mu_{|r|+3} | \; \mathcal{L}_{-1,2,r}[\mathbf{\mu}]^T_{|r|+3}=[\mathbf{0}]^T_{|r|+3} \; \rangle \\
      &=&\langle \; \mu_1 | \; (|r|+2)\mu_1=0 \; \rangle \cong \Z_{|r|+2}.
\end{array}$$
Moreover, using the relations given by $\mathcal{L}_{-1,2,r}$ we
have $\mu_1=\mu_2 \cdots =\mu_{|r|}$ ($\mu_i$'s are the meridians as
before). Therefore, we obtain $$c_1(\xi_{-1,2,r})=PD^{-1}(\mu_1+
\cdots +\mu_{|r|})=PD^{-1}(|r|\mu_1)=|r| \in \Z_{|r|+2}.$$

\medskip \noindent To finish the proof, in each table above we find each
particular case for $(p,q,r)$ such that the corresponding contact
structure $\xi_{p,q,r}$ has binding number 2. Note that the
conditions on $p,q,r$ given in the statement of the theorem excludes
exactly these cases. This completes the proof.
\end{proof}


\section{Remarks on the remaining cases} \label{last section}

\noindent Assume that $r=0,\pm1, |p\,|\geq2, |q|\geq2$. We list all
possible contact structures in Table \ref{r_is_0_-1_1_P_Q}. These
are the only remaining cases from which we still get lens spaces or
their connected sums. Notice that we have already considered the
cases $(-2,q,1)$, and $(2,q,-1)$ in Tables \ref{r_is_1} and
\ref{r_is_-1}, so we do not list them here.

\begin{table}[ht]
  \begin{center}
   \includegraphics{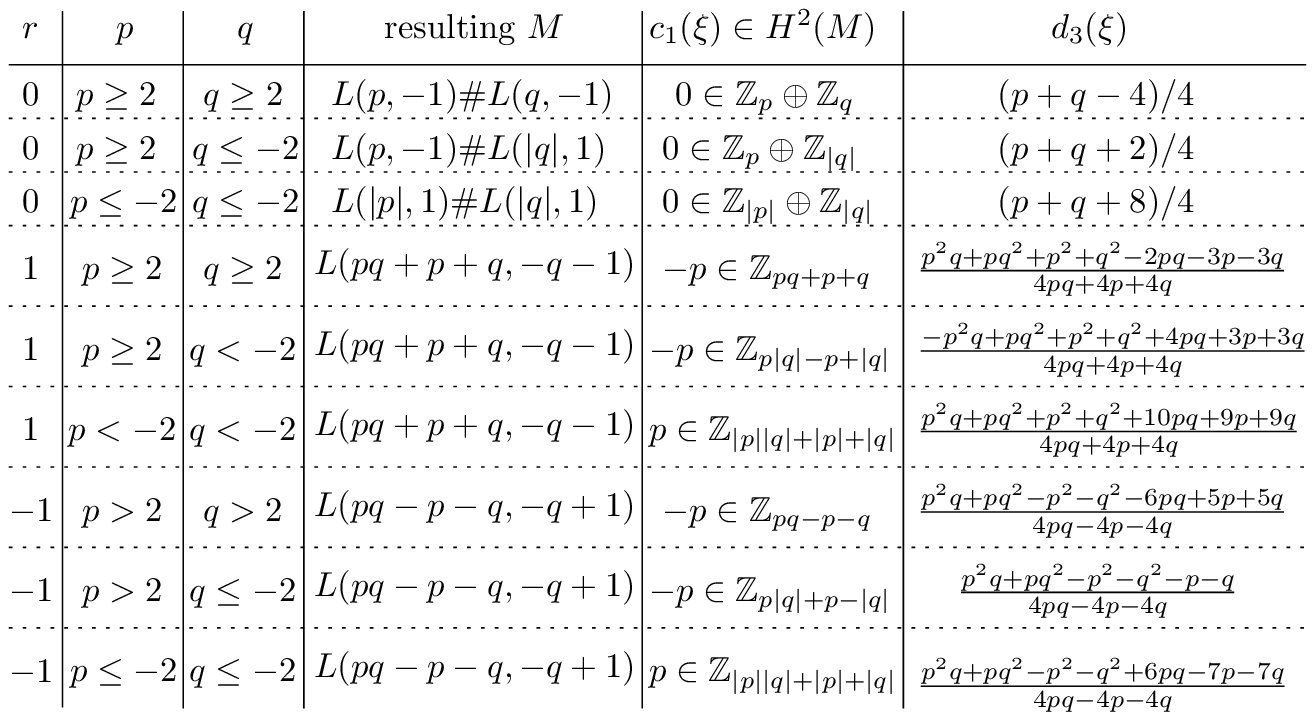}
      \caption{The case $r=0,\pm1$, $|p\,|\geq2$, $|q|\geq2$ ($bn(\xi)=3$ in each row).}
  \label{r_is_0_-1_1_P_Q}
    \end{center}
\end{table}

\noindent As we remarked in Section \ref{introduction} (after
Theorem \ref{the_List}) that one can obtain the complete list
without any repetition: We first simply find all distinct
homeomorphism types of the manifolds which we found in Table
\ref{r_is_0} through Table \ref{r_is_0_-1_1_P_Q}. Then on a fixed
homeomorphism type we compare the pairs $(c_1,d_3)$ coming from the
tables to distinguish the contact structures.

\medskip \noindent  Suppose now that $M$ is a prime Seifert fibered
manifold which is not a lens space. Then as we remarked before we
have $|p\,|\geq2, |q|\geq2$, and $|r|\geq2$. Then two such triples
$(p,q,r), (p',q',r')$ give the same Seifert manifold $Y$ if and only
if
$$\frac{1}{p}+\frac{1}{q}+\frac{1}{r}=\frac{1}{p'}+\frac{1}{q'}+\frac{1}{r'},$$
and $(p',q',r')$ is a permutation of $(p,q,r)$ (see \cite{JN}, for
instance). Notice that we can drop the first condition in our case.
Switching $p$ and $q$ does not change the contact manifold as we
mentioned before. On the other hand, if we switch $r$ and $p$ (or
$r$ and $q$), we might have different contact structures on the same
underlying topological manifold.

\medskip \noindent Another issue is that there are some cases where the first
homology group $H_1(Y(p,q,r))$ is not finite. Indeed, consider the
linking matrix $\mathcal{L}$ of the surgery diagram given on the right in
Figure \ref{our model} as below.
$$\mathcal{L}=\left( \begin{array}{rrrr}
0&1&1&1\\1&p&0&0 \\ 1&0&q&0 \\1&0&0&r \end{array}\right)$$ The
determinant det$(\mathcal{L})=-r(p+q)-pq=0$ implies that $r=-\frac{pq}{p+q}$.
Thus, if $r\neq-\frac{pq}{p+q}$, then $H_1(Y(p,q,r))$ is finite, and
so $d_3(\xi_{p,q,r})$ is still computable since $c_1(\xi_{p,q,r})$
is torsion. For instance, if $p\geq2,q\geq2,r\geq2$ or $p\leq-2,q\leq-2,r\leq-2$, than
det$(\mathcal{L})\neq0$, and so we can distinguish the corresponding
$\xi_{p,q,r}$ by computing the pair $(c_1,d_3)$. Whereas if the
sign of the one of $p,q,r$ is different than the others', then we
might have det$(\mathcal{L})=0$. For instance, for the triples
$(4,4,-2),(3,6,-2)$ and each nonzero integer multiples of them,
det$(A)=0$. So more care is needed for these cases.

\medskip \noindent We would like to end the article by a sample
computation. Assume that det$(\mathcal{L})\neq0$, and that $r\leq2, p\geq2,
q\leq2$ (similar calculations apply for the other cases). We compute
the first homology of $M\approx Y(p,q,r)$ as
$$\begin{array}{lll}
H_1(M)&=&\langle \; \mu_1, \mu_2, \cdots, \mu_{p+q+|r|} | \; \mathcal{L}_{p,q,r}[\mathbf{\mu}]^T_{p+q+|r|}=[\mathbf{0}]^T_{p+q+|r|} \; \rangle \\
      &=&\langle \; \mu_1, \mu_{|r|+1}, \mu_{p+|r|} | \; R_1,R_2,R_3 \;
      \rangle
\end{array}$$
where the relations of the presentation are
$$\begin{array}{cccccccc}
      R_1:& -(2|r|-1)\mu_1&-&(p-1)\mu_{|r|+1}&-&(|q|+1)\mu_{p+|r|}&=&0 \\
      R_2:&              & & p\,\mu_{|r|+1} &-& |q|\,\mu_{p+|r|}  &=&0 \\
      R_3:&    -|r|\,\mu_1&-& p \,\mu_{|r|+1}& &                 &=&0
\end{array}$$
\noindent While getting these relations, we also see that
$\mu_1=\mu_2 \cdots \mu_{|r|}$ (recall $\mu_i$'s are the meridians
to the surgery curves in the family corresponding to $r$ for
$i=1,\cdots, |r|$). Then using this presentation, and knowing that
$c_1(\xi_{p,q,r})=PD^{-1}(|r|\mu_1)$, we can evaluate (understand)
$c_1(\xi_{p,q,r})$ in $H^2(M) \cong H_1(M)$.

\medskip \noindent Now if $c_1(\xi_{p,q,r}) \in H^2(M)$ is a torsion
class, then we can also compute $d_3(\xi_{p,q,r})$ as follows: By
solving the corresponding linear system we get
$$c^2=\frac{p|q||r|}{p|q|+p|r|-|q||r|}.$$ Moreover, we compute
$\sigma(X_{p,q,r})=\sigma(\mathcal{L}_{p,q,r})=-p+|q|+|r|, \chi(X_{p,q,r})=p+|q|+|r|+1$,
and $s=|q|+|r|+1$. Hence, using Corollary \ref{d3_computation}, we
obtain
$$d_3(\xi_{p,q,r})=\frac{8pqr+p^2q+p^2r+4pq^2+4qr^2-pr^2-q^2r-pq-pr-qr}{4pq+4pr+4qr}.$$


\end{document}

%% file: goksty.tex
\def\E{\ifmmode{\mathbb E}\else{$\mathbb E$}\fi} 
\def\N{\ifmmode{\mathbb N}\else{$\mathbb N$}\fi} 
\def\R{\ifmmode{\mathbb R}\else{$\mathbb R$}\fi} 
\def\Q{\ifmmode{\mathbb Q}\else{$\mathbb Q$}\fi} 
\def\C{\ifmmode{\mathbb C}\else{$\mathbb C$}\fi} 
\def\H{\ifmmode{\mathbb H}\else{$\mathbb H$}\fi} 
\def\Z{\ifmmode{\mathbb Z}\else{$\mathbb Z$}\fi} 
\def\P{\ifmmode{\mathbb P}\else{$\mathbb P$}\fi} 
\def\T{\ifmmode{\mathbb T}\else{$\mathbb T$}\fi} 
\def\SS{\ifmmode{\mathbb S}\else{$\mathbb S$}\fi} 
\def\DD{\ifmmode{\mathbb D}\else{$\mathbb D$}\fi} 

\renewcommand{\a}{\alpha}
\renewcommand{\b}{\beta}
\renewcommand{\d}{\delta}
\newcommand{\e}{\varepsilon}
\newcommand{\g}{\gamma}
\newcommand{\G}{\Gamma}
\newcommand{\la}{\lambda}
\newcommand{\La}{\Lambda}
\newcommand{\n}{\nabla}
\newcommand{\var}{\varphi}
\newcommand{\s}{\sigma}
\newcommand{\Sig}{\Sigma}
\renewcommand{\t}{\tau}
\renewcommand{\th}{\theta}
\renewcommand{\O}{\Omega}
\renewcommand{\o}{\omega}
\newcommand{\z}{\zeta}

\newcommand{\ben}{\begin{enumerate}}
\newcommand{\een}{\end{enumerate}}
\newcommand{\be}{\begin{equation}}
\newcommand{\ee}{\end{equation}}
\newcommand{\bea}{\begin{eqnarray}}
\newcommand{\eea}{\end{eqnarray}}
\newcommand{\bc}{\begin{center}}
\newcommand{\ec}{\end{center}}

\newtheorem{thm}{Theorem}[section]
\newtheorem{cor}[thm]{Corollary}
\newtheorem{lem}[thm]{Lemma}
\newtheorem{prop}[thm]{Proposition}
\newtheorem{ax}{Axiom}
\newtheorem{conj}[thm]{Conjecture}

\theoremstyle{definition}
\newtheorem{defn}{Definition}[section]

\theoremstyle{remark}
\newtheorem{rem}{\rm\bfseries{Remark}}[section]
\newtheorem*{notation}{Notation}

\newtheorem{ques}{\rm\bfseries{Question}}[section]
\newtheorem{cons}[rem]{\rm\bfseries{Construction}}
\newtheorem{exm}[rem]{\rm\bfseries{Example}}
